\documentclass{amsart}

\usepackage{amsmath}
\usepackage{amsfonts}
\newtheorem{thm}{Theorem}[section]
\newtheorem{cor}[thm]{Corollary}
\newtheorem{lem}[thm]{Lemma}
\newtheorem{rem}[thm]{Remark}
\newtheorem{prob}[thm]{Problem}
\newtheorem{ex}[thm]{Example}
\newtheorem{prop}[thm]{Proposition}
\numberwithin{equation}{section}
\renewcommand{\Re}{\operatorname{Re}}
\renewcommand{\Im}{\operatorname{Im}}

\title[Carath\'{e}odory interpolation on the non-commutative polydisk]{Carath\'{e}odory interpolation on the non-commutative polydisk}
\author[Dmitry S. Kalyuzhny\u{\i}-Verbovetzki\u{\i}]{Dmitry S. Kalyuzhny\u{\i}-Verbovetzki\u{\i}}
\address{
Department of Mathematics\\
Ben-Gurion University of the Negev\\
P.O. Box 653, Beer-Sheva 84105\\
Israel}
\email{dmitryk@math.bgu.ac.il}
\thanks{The author was supported by the Center for Advanced Studies in Mathematics, Ben-Gurion University of the Negev.}

\begin{document}
\begin{abstract}
The Carath\'{e}odory problem in the $N$-variable non-commutative Herglotz--Agler class and the Carath\'{e}odory--Fej\'{e}r problem in the $N$-variable non-commutative Schur--Agler class are posed. It is shown that the Carath\'{e}o\-dory (resp., Carath\'{e}odory--Fej\'{e}r) problem has a solution if and only if the non-commutative polynomial with given operator coefficients (the data of the problem indexed by an admissible set $\Lambda$) takes operator values with positive semidefinite real part (resp., contractive operator values) on $N$-tuples of $\Lambda$-jointly nilpotent contractive $n\times n$ matrices, for all $n\in\mathbb{N}$. 
\end{abstract}
\maketitle

\section{Introduction}\label{sec:intro}
The classical Carath\'{e}odory interpolation problem is the following: given a sequence of complex numbers $c_0>0,\ c_1,\ldots,c_m$, find a holomorphic function $f(z)=f_0+f_1z+f_2z^2+\cdots$ on the open unit disk $\mathbb{D}$ whose values in $\mathbb{D}$ have positive real part (i.e., $f$ belongs to the \emph{Herglotz}, or \emph{Carath\'{e}odory, class} $\mathcal{H}_1$, where the subscript $1$ stands for the one-variable case) such that 
$$f_0=\frac{c_0}{2},\quad f_1=c_1, \ldots,  f_m=c_m.$$
This problem has been posed by Constantin Carath\'{e}odory in \cite{Car1,Car2} where the criteria of its solvability and of the uniqueness of its solution were presented. Toeplitz has noticed in \cite{Toep} that the original solvability criterion from \cite{Car1}, which was formulated in terms of convex bodies, admits the following formulation in terms of the coefficients $c_k,\ k=0,\ldots,m$:
the Carath\'{e}odory problem for these data has a solution if and only if the $(m+1)\times (m+1)$ matrix
	\begin{equation}\label{t}
	T_c=\left[
\begin{array}{lllll}
	c_0       & c_1^*    & \ldots & c_{m-1}^* & c_m^*\\
	c_1     & \ddots    & \ddots    & \ldots  & c_{m-1}^*\\
	\vdots    & \ddots & \ddots & \ddots  & \vdots\\
		c_{m-1} & \ldots & \ddots & \ddots     & c_1^*\\
	c_m     & c_{m-1} & \ldots & c_1 & c_0  
\end{array}\right] 
\end{equation}
 is positive semidefinite (here $c_k^*=\overline{c_k}$). From the integral representation by Riesz \cite{Ri} and Herglotz \cite{He}
 \begin{equation}\label{rh}
 f(z)=\frac{1}{2}\int_\mathbb{T}\frac{1+\bar{\lambda}z}{1-\bar{\lambda}z}\,d\mu(\lambda)+i\Im f(0),\quad z\in\mathbb{D},
 \end{equation}
 which characterizes functions from $\mathcal{H}_1$ (here $\mu$ is a positive Borel measure on the unit circle $\mathbb{T}$; in the case where $f(0)=\frac{1}{2}$ the second term in the right-hand side of \eqref{rh} is dropped out and $\mu$ has full variation $|\mu|=1$) one obtains a representation for the Taylor coefficients of $f\in\mathcal{H}_1$:
 $$f_0=\frac{|\mu|}{2}+i\Im f(0),\quad f_k=\int_\mathbb{T}\bar{\lambda}^k\,d\mu(\lambda),\quad k=1,2,\ldots .$$
Thus, the Carath\'{e}odory problem has a solution if and only if there exists a positive Borel measure $\mu$ on $\mathbb{T}$ such that
\begin{equation}\label{mom}
c_k=\int_\mathbb{T}\bar{\lambda}^k\,d\mu(\lambda),\quad k=0,\ldots, m,
\end{equation}
 i.e., $\mu$ solves the \emph{trigonometric moment problem} for the data $c_k,\ k=0,\ldots,m$.
 
 In the operator case the data of the Carath\'{e}odory problem are bounded linear operators $c_0\geq 0,\ c_1,\ldots,c_m$ on a separable Hilbert space\footnote{In this paper we will consider separable Hilbert spaces only, and omit ``separable" for brevity.} $\mathcal{Y}$, and the class $\mathcal{H}_1$ is replaced by the class $\mathcal{H}_1(\mathcal{Y})$ of holomorphic functions on $\mathbb{D}$ whose values are bounded linear operators\footnote{For Hilbert spaces $\mathcal{Y}$ and $\mathcal{H}$, we shall use the notation $\mathcal{L(Y,H)}$ (resp., $\mathcal{L(Y)}$) for the Banach space of bounded linear operators from $\mathcal{Y}$ to $\mathcal{H}$ (resp., from $\mathcal{Y}$ to itself).} on $\mathcal{Y}$ with positive semidefinite real part. Then the Carath\'{e}odory--Toeplitz criterion, representation \eqref{rh} for $f\in\mathcal{H}_1(\mathcal{Y})$, and trigonometric moment representation \eqref{mom} hold true with the operator block matrix $T_c$ in \eqref{t}, a positive Borel $\mathcal{L(Y)}$-valued measure $\mu$, and the convergence of integrals in \eqref{rh} and \eqref{mom} in the strong operator topology. Riesz--Herglotz representation \eqref{rh} for the case where $f(0)=\frac{I_\mathcal{Y}}{2}$, and thus moment representation \eqref{mom} for the case where $c_0=I_\mathcal{Y}$ admit the following operator form:
\begin{eqnarray}
f(z) &=& \frac{1}{2}V^*(I_\mathcal{H}+zG)(I_\mathcal{H}-zG)^{-1}V,\quad z\in\mathbb{D}, \label{rh-o}\\
c_k  &=& V^*G^kV,\qquad k=0,\ldots,m,\label{mom-o}
\end{eqnarray}
 where $G$ is a unitary operator on some auxiliary Hilbert space $\mathcal{H}$, and $V\in\mathcal{L(Y,H)}$ is an isometry. These results are due to Neumark \cite{Ne}.
 
A similar problem was considered first by Carath\'{e}odory and Fej\'{e}r in \cite{CF} for the \emph{Schur class} $\mathcal{S}_1$ of holomorphic contractive functions on $\mathbb{D}$ in the place of the Herglotz class $\mathcal{H}_1$:  given a sequence of complex numbers $s_0,\ldots,s_m$, find a holomorphic function $F(z)=F_0+F_1z+F_2z^2+\cdots$ from the class $\mathcal{S}_1$ such that 
$$F_0=s_0,\ldots,F_m=s_m.$$
Schur has proved in \cite{Sch} that the Carath\'{e}odory--Fej\'{e}r problem has a solution if and only if the matrix
	\begin{equation}\label{st}
	T_s=\left[
\begin{array}{llll}
	s_0       & 0    & \ldots & 0\\
	s_1     & \ddots   & \ddots    & \vdots \\
	\vdots    & \ddots & \ddots & 0\\
	s_m     & \ldots  & s_1 & s_0  
\end{array}\right] 
\end{equation}
is contractive\footnote{Matrix norm considered in this paper is operator $(2,2)$-norm, i.e., the maximal singular value of a matrix.}, i.e., $\| T_s\|\leq 1$. In the operator case the data of the Carath\'{e}odory--Fej\'{e}r problem are operators $s_0,\ldots,s_m\in\mathcal{L(U,Y)}$, with  Hilbert spaces $\mathcal{U}$ and $\mathcal{Y}$, the class $\mathcal{S}_1$ is replaced by the class $\mathcal{S}_1(\mathcal{U,Y})$ of holomorphic functions on $\mathbb{D}$ whose values are contractive operators from $\mathcal{L(U,Y)}$ and the Schur criterion is formulated in the same way as in the scalar case, with the operator block matrix $T_s$.

Let us note that a common operatorial view at the Carath\'{e}odory problem, Carath\'{e}odory--Fej\'{e}r problem, and their relative Nevanlinna--Pick problem, and a certain operator dilation scheme which unifies these problems were first presented in the fundamental paper of Sarason \cite{Sar}. These ideas in an abstract form have been expressed in the commutant lifting theorem of Sz.-Nagy and Foia\c{s} (see \cite{SzNF}) which is used now as one of the approaches to various interpolation problems. For further details on the classical and operator versions of the Carath\'{e}odory problem and other interpolation problems, see \cite{Akh, KN, FF,Co, Sakh}.

There exist various generalizations of the Carath\'{e}odory problem and other interpolation problems to the case of several complex variables, depending on the type of a classical domain in $\mathbb{C}^N$ serving as a counterpart of $\mathbb{D}$ and on the class of interpolating functions. Due to a version of the Riesz--Herglotz formula \eqref{rh} for
the unit polydisk $\mathbb{D}^N$ obtained by Kor\'{a}nyi and Puk\'{a}nszky \cite{KP}, one can characterize the coefficients of a function from the multivariable \emph{Herglotz class} $\mathcal{H}_N(\mathcal{Y})$ (the class of holomorphic functions on $\mathbb{D}^N$ taking operator values from $\mathcal{L(Y)}$ with positive semidefinite real part) in terms of a $\mathcal{L(Y)}$-valued positive Borel measure $\mu$ whose Fourier coefficients with multi-indices outside $\mathbb{Z}^N_+$ and $\mathbb{Z}^N_-$, the positive and the negative discrete octants, are zero. However, an appropriate multivariable analogue of \eqref{rh-o} (and thus, of \eqref{mom-o}) can be obtained either for the case $N=2$ or for the subclass $\mathcal{HA}_N(\mathcal{Y})\subset\mathcal{H}_N(\mathcal{Y})$ which is proper for $N>2$. The latter subclass, which is called the  \emph{Herglotz--Agler class}, has been introduced by Agler in \cite{Ag1}, where the analogue of \eqref{rh-o} has been obtained. This class $\mathcal{HA}_N(\mathcal{Y})$ consists of holomorphic $\mathcal{L(Y)}$-valued functions on $\mathbb{D}^N$ whose values on any $N$-tuple of commuting strict contractions on a common Hilbert space (in the sense of hereditary functional calculus introduced in \cite{Ag1}) have positive semidefinite real part. Some partial results on the Carath\'{e}odory--Fej\'{e}r problem in the class $\mathcal{S}_N$ (the \emph{Schur class} of contractive holomorphic functions on $\mathbb{D}^N$) have been obtained in \cite{Pf,DKh}. The Carath\'{e}odory and Carath\'{e}odory--Fej\'{e}r problems in the Herglotz--Agler class $\mathcal{HA}_N(\mathcal{Y})$
and the \emph{Schur--Agler class $\mathcal{SA}_N(\mathcal{U,Y})$} (the class of holomorphic $\mathcal{L(U,Y)}$-valued functions on $\mathbb{D}^N$ which take contractive operator values on any $N$-tuple of commuting strict contractions, in the sense of Agler's hereditary functional calculus), respectively, were studied in \cite{EPP,BLTT,W}. Various versions of the Sarason theorem, the Sz.-Nagy--Foia\c{s} commutant lifting theorem and the Nevanlinna--Pick problem on $\mathbb{D}^N$, in the classes $\mathcal{H}_N(\mathcal{Y}), \mathcal{HA}_N(\mathcal{Y}), \mathcal{S}_N(\mathcal{U,Y})$ and $\mathcal{SA}_N(\mathcal{U,Y})$ were considered in \cite{CS,CW,AgM1,BT,BLTT,EPP}. The Kor\'{a}nyi--Puk\'{a}nszky version of the Riesz--Herglotz formula has been generalized in \cite{AD} to a wide class of domains in $\mathbb{C}^N$ which contains, in particular, all classical symmetric domains. Certain partial results on the Carath\'{e}odory--Fej\'{e}r problem for bounded full circular domains in $\mathbb{C}^N$ can be found in \cite{DKh}. A generalization of  
the Agler representation theorem from \cite{Ag1} to a class of so-called polynomially defined domains in $\mathbb{C}^N$ has been obtained in \cite{AT,BB1} where also the Nevanlinna--Pick problem in the Schur--Agler class of functions on such a domain was studied. The Nevanlinna--Pick and Carath\'{e}odory--Fej\'{e}r problems in the class of contractive multipliers on the reproducing kernel Hilbert space of holomorphic functions on the \emph{unit ball} $\mathbb{B}_N:=\{ z\in\mathbb{C}^N:\ \sum_{k=1}^N|z_k|^2<1\}$, with the reproducing kernel $k_N(z,z')=\frac{1}{1-\left\langle z,z'\right\rangle}$, or more generally, on the reproducing kernel Hilbert space of functions on a set $\Omega$, with the reproducing kernel whose reciprocal has exactly one positive square, were studied starting with the unpublished paper of Agler \cite{Ag2} by many authors (e.g., see \cite{Ag3,McC1,McC2,Q,AgM2,BTV,BB2}). Let us mention also the approach to interpolation problems on $\mathbb{B}_N$ via the commutant lifting theorem in the non-commutative setting of the Toeplitz algebra of operators acting on the Fock space by Popescu \cite{P1,P2} and subsequent use of symmetrization argument (see \cite{P3,AP,DP}). In this non-commutative setting the Carath\'{e}odory--Fej\'{e}r problem was studied in \cite{P3,P4,CoJ}. (A certain generalization of Popescu's non-commutative setting and a more general Nevanlinna--Pick interpolation problem appears in a recent paper \cite{MS}.) Let us remark that one can interprete the latter results in terms of functions on the \emph{non-commutative unit ball} $\mathcal{B}_N$ which is the collection of \emph{strict row contractions}, i.e., $N$-tuples of bounded linear operators $\mathbf{T}=(T_1,\ldots,T_N)$ on a common Hilbert space $\mathcal{E}$ such that $\sum_{k=1}^NT_kT_k^*<I_\mathcal{E}$.

In the present paper, we are working on another domain, the \emph{non-commutative unit polydisk} $\mathcal{D}^N$ which is the collection of $N$-tuples $\mathbf{T}=(T_1,\ldots,T_N)$ of strict contractions on a common Hilbert space $\mathcal{E}$, i.e, $\| T_k\|<1,\ k=1,\ldots,N$, or on the \emph{non-commutative matrix unit polydisk}, which is a subbdomain $\mathcal{D}^N_{\rm matr}\subset\mathcal{D}^N$ consisting of $N$-tuples of strict contractions on $\mathbb{C}^n$, for all $n\in\mathbb{N}$. The domain $\mathcal{D}^N$ is a special case of a bit more general non-commutative domain $\mathcal{D}_G$ considered in the recent paper of Ball, Groenewald and Malakorn \cite{BGM} where the \emph{non-commutative Schur--Agler class} $\mathcal{SA}_N^{\rm nc}(\mathcal{U,Y})$ was introduced and studied in the framework of structured non-commutative multidimensional conservative linear systems. The domain $\mathcal{D}^N_{\rm matr}$ appears in \cite{AK2}. We consider non-commutative formal power series which converge on $\mathcal{D}^N$. We introduce the \emph{non-commutative Herglotz--Agler class} $\mathcal{HA}_N^{\rm nc}(\mathcal{Y})$ of such series which take on $\mathcal{D}^N$ operator values with positive semidefinite real part, and study the Carath\'{e}odory problem in this class, as well as the Carath\'{e}odory--Fej\'{e}r problem in the class $\mathcal{SA}_N^{\rm nc}(\mathcal{U,Y})$.

To give an idea on our main results, criteria of solvability of these problems, let us first come back to the one-variable case. Let $S$ denote the standard shift $(m+1)\times(m+1)$ matrix:
\begin{equation}\label{s}
S=\left[\begin{array}{cccc}
0 &        &     & \\
1 & \ddots &     & \\
  & \ddots & \ddots & \\
  &        &     1  & 0 
\end{array}\right],
\end{equation}
i.e., $S_{ij}=1$ for $i-j=1$, and $S_{ij}=0$ otherwise. Then
$$T_c=I_{m+1}\otimes c_0+\sum_{k=1}^mS^k\otimes c_k+\sum_{k=1}^mS^{*k}\otimes c_k^*.$$
If one defines
$$p(z):=\frac{c_0}{2}+\sum_{k=1}^mc_kz^k,\qquad z\in\mathbb{C},$$
then the Carath\'{e}odory--Toeplitz criterion can be formulated as the positive semidefiniteness of the operator $2\Re p^{\rm l}(S)=p^{\rm l}(S)+p^{\rm l}(S)^*$, where
$$p^{\rm l}(S):=\frac{I_{m+1}\otimes c_0}{2}+\sum_{k=1}^mS^k\otimes c_k,$$
or equivalently, of the operator $2\Re p(S)$, where
$$p(S)=p^{\rm r}(S):=\frac{c_0\otimes I_{m+1}}{2}+\sum_{k=1}^mc_k\otimes S^k$$
(we shall usually omit the superscript ``r", however keep the superscript ``l" when we use the writing of a polynomial with powers on the left). By \cite[Section~2.5]{Arv2}, any contraction $T$ on a Hilbert space $\mathcal{E}$ which is \emph{nilpotent of rank at most $m+1$}, i.e., such that 
$$T^k=0,\qquad k=m+1,m+2,\ldots,$$
admits a \emph{dilation} of the form $S\otimes I_\mathcal{H}$, with some Hilbert space $\mathcal{H}$, i.e., there exists an isometry $V\in\mathcal{L(E},\mathbb{C}^{m+1}\otimes\mathcal{H})$ such that
$$T^k=V^*(S^k\otimes I_\mathcal{H})V,\qquad k=1,2,\ldots .$$
Since $S$ is a nilpotent matrix with rank of nilpotency $m+1$, we obtain the following criterion:
\emph{the Carath\'{e}odory problem with data $c_0\geq 0,c_1,\ldots,c_m\in\mathcal{L(Y)}$ has a solution if and only if $\Re p(T)\geq 0$ for every nilpotent operator $T$ with rank of nilpotency at most $m+1$.} 
Analogously, $$T_s=\sum_{k=0}^mS^k\otimes s_k.$$
If one defines
$$q(z):=\sum_{k=0}^ms_kz^k,\qquad z\in\mathbb{C},$$
then the Schur criterion can be formulated as the contractivity of the operator $q^{\rm l}(S)$ where
$$q^{\rm l}(S)=\sum_{k=0}^mS^k\otimes s_k,$$
or equivalently, of the operator $q(S)$, where
$$q(S)=\sum_{k=0}^ms_k\otimes S^k.$$
Thus, \emph{the Carath\'{e}odory--Fej\'{e}r problem with data $s_0,\ldots,s_m\in\mathcal{L(U,Y)}$ has a solution if and only if $\| q(T)\|\leq 1$ for every nilpotent operator $T$ with rank of nilpotency at most $m+1$.}

The main results of the paper are generalizations of these criteria to the multivariable non-commutative case, where the positivity or contractivity of a non-commutative polynomial is tested on $N$-tuples of \emph{jointly nilpotent} contractions $\mathbf{T}=(T_1,\ldots,T_N)$, i.e, $\| T_k\|\leq 1,\  k=1,\ldots,N,$ and $T_{i_1}\cdots T_{i_k}=0$ outside some finite set of strings $(i_1,\ldots,i_k),\ k\in\mathbb{N}$. To obtain these criteria, we first deduce the analogue of \eqref{rh-o} for non-commutative formal power series of the class $\mathcal{HA}_N^{\rm nc}(\mathcal{Y})$ from the realization formula obtained in \cite{BGM} for the class $\mathcal{SA}_N^{\rm nc}(\mathcal{Y})$. Thus, an anologue of \eqref{mom-o} is also obtained. A counterpart of unitary operator $G$ from \eqref{rh-o} and \eqref{mom-o} is $N$-tuple
$\mathbf{G}=(G_1,\ldots,G_N)$ of bounded linear operators on a common Hilbert space satisfying the following condition:
\begin{equation}\label{g}
\zeta\mathbf{G}:=\sum_{k=1}^N\zeta_kG_k\quad {\rm is\ unitary\ for\ every}\ \zeta\in\mathbb{T}^N.
\end{equation} 
We denote by $\mathcal{G}_N$ the class of such $N$-tuples of operators. Note that an $N$-tuple $\mathbf{G}$ from the class  $\mathcal{G}_N$ appears also in Agler's representation formula for functions from the (commutative) class $\mathcal{HA}_N(\mathcal{Y})$ in \cite{Ag1}, and in a realization formula for functions from the subclass $\mathcal{SA}_N^0(\mathcal{U,Y})\subset\mathcal{SA}_N(\mathcal{U,Y})$ which consists of functions vanishing at zero, in \cite{K2}. We deduce the criterion of solvability of the Carath\'{e}odory--Fej\'{e}r problem in the class $\mathcal{SA}_N^{\rm nc}(\mathcal{U,Y})$  from the one for the Carath\'{e}odory problem in the class $\mathcal{HA}_N^{\rm nc}(\mathcal{Y})$, which we obtain first.

The main tools in our work, besides the realization formula from \cite{BGM} mentioned above, are the following: the properties of operator $N$-tuples from the class $\mathcal{G}_N$ which have been established in \cite{K2} and some their new properies which we obtain in the present paper; the factorization result of McCullough \cite{McC3} for non-commutative hereditary polynomials; the Arveson extension theorem \cite{Arv1}; the Stinespring representation theorem \cite{St} for completely positive maps of $C^*$-algebras; the Sz.-Nagy and Foia\c{s} theorem on the existence of a unitary dilation of an $N$-tuple of (not necessarily commuting) contractions \cite{SzNF}; the Amitsur--Levitzki theorem  on the non-existence of non-commutative polynomial relations valid for infinitely many matrix rings $\mathbb{C}^{n_j\times n_j},\ j=1,2,\ldots$ (see, e.g., \cite[pp. 22--23]{Ro}).

The structure of the paper is the following. In Section~\ref{sec:tuples} we study certain classes of operator $N$-tuples. In particular, we establish duality properties of the classes $\mathcal{G}_N$ and $\mathcal{U}^N$. The latter is the class of $N$-tuples $\mathbf{U}=(U_1,\ldots,U_N)$ of unitary operators on a common Hilbert space, which serves as another generalization (in addition to $\mathcal{G}_N$) of the class of single unitaries. This duality is observed also in Lemma~\ref{lem:main} which is proved in Section~\ref{sec:prob}. In Section~\ref{sec:ha} we introduce and characterize the non-commutative Herglotz--Agler class $\mathcal{HA}_N^{\rm nc}(\mathcal{Y})$. In Section~\ref{sec:prob} we formulate the Carath\'{e}odory problem in the class  $\mathcal{HA}_N^{\rm nc}(\mathcal{Y})$ and prove the necessary and sufficient conditions for its solvability. In Section~\ref{sec:prob-cf} we formulate the Carath\'{e}odory--Fej\'{e}r problem in the class  $\mathcal{SA}_N^{\rm nc}(\mathcal{U,Y})$ and obtain a criterion of its solvability.
 
\section{Some classes of operator $N$-tuples}\label{sec:tuples}
Let us define the classes of $N$-tuples of operators which are considered in this paper. In addition to the classes $\mathcal{D}^N$, $\mathcal{D}^N_{\rm matr}$, $\mathcal{G}_N$ and $\mathcal{U}^N$ already mentioned in Section~\ref{sec:intro}, let us define the class $\mathcal{C}^N$ which consists of $N$-tuples $\mathbf{C}=(C_1,\ldots,C_N)$ of contractions on a common Hilbert space, i.e., $\| C_k\|\leq 1,\ k=1,\ldots,N$. The class $\mathcal{G}_N$ is characterized by the following proposition which is a consequence of \cite[Proposition 2.4]{K2}.
\begin{prop}\label{prop:gn}
Let $\mathbf{G}=(G_1,\ldots,G_N)\in\mathcal{L(H)}^N$, with a Hilbert space $\mathcal{H}$. The following statements are equivalent:
\allowdisplaybreaks
\begin{description}
\item[(i)] $\mathbf{G}\in\mathcal{G}_N$;
	\item[(ii)] $\mathbf{G}$ satisfies the conditions
	\begin{eqnarray}
	\sum_{k=1}^NG_k^*G_k &=& I_\mathcal{H},\label{i1}\\
	G_k^*G_j &=& 0,\qquad k\neq j,\label{i2}\\
	\sum_{k=1}^NG_kG_k^* &=& I_\mathcal{H},\label{i3}\\
	G_kG_j^* &=& 0,\qquad k\neq j;\label{i4}
		\end{eqnarray} 
		\item[(iii)] the operator $G^0:=\sum_{k=1}^NG_k$ is unitary, and there exists a resolution of identity $I_\mathcal{H}=\sum_{k=1}^NP_k^-$, where $(P_k^-)^2=P_k^-=(P_k^-)^*,\ k=1,\ldots,N$, and $P_k^-P_j^-=0$ for $k\neq j$, such that
		$$G_k=G^0P_k^-,\qquad k=1,\ldots,N;$$
		 	\item[(iv)] the operator $G^0:=\sum_{k=1}^NG_k$ is unitary, and there exists a resolution of identity $I_\mathcal{H}=\sum_{k=1}^NP_k^+$, where $(P_k^+)^2=P_k^+=(P_k^+)^*,\ k=1,\ldots,N$, and $P_k^+P_j^+=0$ for $k\neq j$, such that
		$$G_k=P_k^+G^0,\qquad k=1,\ldots,N;$$
				\item[(v)] for every $\mathbf{U}=(U_1,\ldots,U_N)\in\mathcal{U}^N\cap\mathcal{L}(\mathcal{K})^N$, with a Hilbert space $\mathcal{K}$, the operator
		\begin{equation}\label{l-pencil}
		\mathbf{U\otimes G}:=\sum_{k=1}^NU_k\otimes G_k\in\mathcal{L}(\mathcal{K}\otimes\mathcal{H})
		\end{equation}
		is unitary;
		\item[(vi)] for every $\mathbf{U}=(U_1,\ldots,U_N)\in\mathcal{U}^N\cap\mathcal{L}(\mathcal{K})^N$, with a Hilbert space $\mathcal{K}$, the operator
		\begin{equation}\label{r-pencil}
		\mathbf{G\otimes U}:=\sum_{k=1}^NG_k\otimes U_k\in\mathcal{L}(\mathcal{H\otimes K})
		\end{equation}
		is unitary.
		\end{description}
\end{prop}
\begin{cor}\label{cor:gn}
If $\mathbf{G}\in\mathcal{G}_N\cap\mathcal{L(H)}^N$ then: 
\begin{description}
\item[(a)] for every $\mathbf{C}\in\mathcal{C}^N\cap\mathcal{L}(\mathcal{E})^N$, with a Hilbert space $\mathcal{E}$, the operators $\mathbf{G\otimes C}$ and $\mathbf{C\otimes G}$ are contractions; 
\item[(b)]  for every $\mathbf{C}\in\mathcal{D}^N\cap\mathcal{L}(\mathcal{E})^N$, with a Hilbert space $\mathcal{E}$, the operators $\mathbf{G\otimes C}$ and $\mathbf{C\otimes G}$ are strict contractions.
\end{description}
\end{cor}
\begin{proof}
(a) Any $\mathbf{C}\in\mathcal{C}^N\cap\mathcal{L}(\mathcal{E})^N$ has a \emph{unitary dilation} \cite{SzNF}, i.e., there exists an $N$-tuple  $\mathbf{U}\in\mathcal{U}^N\cap\mathcal{L(K)}^N$, with some Hilbert space $\mathcal{K}\supset\mathcal{E}$, such that
$$C_{i_1}\cdots C_{i_l}=P_\mathcal{E}U_{i_1}\cdots U_{i_l}\big|_\mathcal{E}, \qquad l\in\mathbb{N},\ i_1,\ldots,i_l\in\{ 1,\ldots,N\},
$$
where $P_\mathcal{E}$ denotes the orthogonal projection onto the subspace $\mathcal{E}$ in $\mathcal{K}$. Therefore, $$\mathbf{G\otimes C}=\sum_{k=1}^NG_k\otimes P_\mathcal{E}U_k\big|_\mathcal{E}=(I_\mathcal{H}\otimes P_\mathcal{E})(\mathbf{G\otimes U})\big|_{\mathcal{H}\otimes\mathcal{E}},$$ and since by Proposition~\ref{prop:gn} $\mathbf{G\otimes U}$ is unitary, $\mathbf{G\otimes C}$ is a contraction. Analogously, $\mathbf{C\otimes G}$ is a contraction.

(b) If $\mathbf{C}\in\mathcal{D}^N\cap\mathcal{L}(\mathcal{E})^N$ is non-zero (otherwise the statement is trivial) then $\widetilde{\mathbf{C}}:=(\max_{1\leq k\leq N}\| C_k\| )^{-1}\mathbf{C}\in\mathcal{C}^N\cap\mathcal{L}(\mathcal{E})^N$. By part (a) of this Proposition, $\mathbf{G}\otimes \widetilde{\mathbf{C}}$ and $\widetilde{\mathbf{C}}\otimes \mathbf{G}$ are contractions. Therefore, 
 $\mathbf{G\otimes C}$ and $\mathbf{C\otimes G}$ are strict contractions with norm at most $\max_{1\leq k\leq N}\| C_k\|$.
\end{proof}
The following proposition is dual to Proposition~\ref{prop:gn}.
\begin{prop}\label{prop:un} 
Let $\mathbf{U}=(U_1,\ldots,U_N)\in\mathcal{L(K)}^N$, with a Hilbert space $\mathcal{K}$. The following statements are equivalent:
\begin{description}
	\item[(i)] $\mathbf{U}\in\mathcal{U}^N$;
	\item[(ii)] for every $\mathbf{G}\in\mathcal{G}_N\cap\mathcal{L}(\mathcal{H})^N$, with a Hilbert space $\mathcal{H}$, the operator $\mathbf{U\otimes G}\in\mathcal{L}(\mathcal{K}\otimes\mathcal{H})^N$ is unitary;
	\item[(iii)] for every $\mathbf{G}\in\mathcal{G}_N\cap\mathcal{L}(\mathcal{H})^N$, with a Hilbert space $\mathcal{H}$, the operator $\mathbf{G\otimes U}\in\mathcal{L}(\mathcal{H}\otimes\mathcal{K})^N$ is unitary.
\end{description}
\end{prop}
\begin{proof}
If $\mathbf{U}\in\mathcal{U}^N$ and $\mathbf{G}\in\mathcal{G}_N$ then implications (i)$\Rightarrow$(ii) and (i)$\Rightarrow$(iii) follow from Proposition~\ref{prop:gn}.
For the proof of (ii)$\Rightarrow$(i) and (iii)$\Rightarrow$(i), one can choose $$\mathbf{G}^{(k)}:=(0,\ldots,0,1,0,\ldots, 0)\in\mathbb{C}^N\cong\mathcal{L}(\mathbb{C})^N,$$ where $1$ is on the $k$'th position and $0$ is on the other positions. It is clear that $ \mathbf{G}^{(k)}\in\mathcal{G}_N$. Since for every $k\in\{ 1,\ldots,N\}$ the operator $U_k=\mathbf{U\otimes G}^{(k)}=\mathbf{G}^{(k)}\otimes\mathbf{U}$ is unitary, $\mathbf{U}\in\mathcal{U}^N$. 
\end{proof}
For $N$-tuples of operators $\mathbf{X}=(X_1,\ldots,X_N)\in\mathcal{L(X)}^N$ and $\mathbf{Y}=(Y_1,\ldots,Y_N)\in\mathcal{L(Y)}^N$ on Hilbert spaces $\mathcal{X}$ and $\mathcal{Y}$, respectively, define their \emph{Schur tensor product} as the $N$-tuple of operators 
\begin{equation}\label{schur}
\mathbf{X}\stackrel{\circ}{\otimes}\mathbf{Y}:=(X_1\otimes Y_1,\ldots,X_N\otimes Y_N)\in\mathcal{L(X\otimes Y)}^N.
\end{equation}
\begin{prop}\label{prop:schur}
For any $\mathbf{G}\in\mathcal{G}^N$ and $\mathbf{U}\in\mathcal{U}^N$ both $\mathbf{G}\stackrel{\circ}{\otimes}\mathbf{U}$ and $\mathbf{U}\stackrel{\circ}{\otimes}\mathbf{G}$ belong to the class $\mathcal{G}^N$.
\end{prop}
\begin{proof}
For an arbitrary $\tilde{\mathbf{U}}\in\mathcal{U}^N$ the operator
$$(\mathbf{G}\stackrel{\circ}{\otimes}\mathbf{U})\otimes\tilde{\mathbf{U}}=\sum_{k=1}^NG_k\otimes U_k\otimes \tilde{U}_k=\mathbf{G}\otimes (\mathbf{U}\stackrel{\circ}{\otimes}\tilde{\mathbf{U}})$$
is unitary, by Proposition~\ref{prop:gn} and due to the fact that $U_k\otimes \tilde{U}_k$ are unitary operators for all $k=1,\ldots,N$, i.e., $\mathbf{U}\stackrel{\circ}{\otimes}\tilde{\mathbf{U}}\in\mathcal{U}^N$. Thus, again by Proposition~\ref{prop:gn}, $\mathbf{G}\stackrel{\circ}{\otimes}\mathbf{U}\in\mathcal{G}^N$. Analogously, $\mathbf{U}\stackrel{\circ}{\otimes}\mathbf{G}\in\mathcal{G}^N$. 
\end{proof}
Let us note that for the classes introduced above the following inclusions hold: 
$$\begin{array}{ccccccc}
                         &           &               &         & \mathcal{G}_N &         &       \\        
                         &           &               &         &       \cap    &         &       \\  
\mathcal{D}^N_{\rm matr} & \subset   & \mathcal{D}^N & \subset & \mathcal{C}^N & \supset & \mathcal{U}^N.
\end{array}$$
A couple of additional classes of operator $N$-tuples will be considered in Section~\ref{sec:prob}.

\section{The non-commutative Herglotz--Agler class}\label{sec:ha}
 Let us give some necessary definitions in our non-commutative setting. The \emph{free semigroup $\mathcal{F}_N$ with the generators $g_1,\ldots,g_N$ (the letters)} and the neutral element $\emptyset$ (\emph{empty word}) has a product defined as follow: if two its elements (\emph{words}) are given by $w=g_{i_1}\cdots g_{i_m}$ and  $w'=g_{j_1}\cdots g_{j_{m'}}$ then their product is $ww'=g_{i_1}\cdots g_{i_m}g_{j_1}\cdots g_{j_{m'}},$
and  $w\emptyset=\emptyset w=w$. The \emph{length of the word} $w=g_{i_1}\cdots g_{i_m}$ is $|w|=m$, and $|\emptyset|=0$. 
The non-commutative algebra $\mathcal{L(Y)}\left\langle \left\langle z_1,\ldots,z_N\right\rangle\right\rangle$ consists 
of \emph{formal power series} $f$ with coefficients $f_w\in\mathcal{L(Y)},\ w\in\mathcal{F}_N$, for a Hilbert space $\mathcal{Y}$, of the form
$$f(z)=\sum_{w\in\mathcal{F}_N}f_wz^w,$$
where for the indeterminates $z=(z_1,\ldots,z_N)$ and  words $w=g_{i_1}\cdots g_{i_m}$ one sets $z^w=z_{i_1}\cdots z_{i_m},\ z^\emptyset =1$. We assume that indeterminates $z_k$ formally commute with coefficients $f_w$. A formal power series $f$ is invertible in this algebra if and only if $f_\emptyset$ is inverible. Indeed, if $f(z)\phi(z)=\phi(z)f(z)=I_\mathcal{Y}$ then $f_\emptyset\phi_\emptyset=\phi_\emptyset f_\emptyset=I_\mathcal{Y}$, i.e., $\phi_\emptyset=f_\emptyset^{-1}$. Conversely, if $f_\emptyset$ is inverible then the series
 $$
\phi(z)=\sum_{k=0}^\infty \left(I_\mathcal{Y}-
f_\emptyset^{-1}f(z)\right)^kf_\emptyset^{-1}.
$$
is the inverse of $f$.
This formal power series is well defined since the
expansion of $\left(I_\mathcal{Y}-f_\emptyset^{-1}f\right)^k$ contains words of length
at least $k$, and thus the expressions for coefficients  $\phi_w$ are finite
sums. The subalgebra $\mathcal{L(Y)}\left\langle z_1,\ldots,z_N\right\rangle$ of the algebra $\mathcal{L(Y)}\left\langle \left\langle z_1,\ldots,z_N\right\rangle\right\rangle$ consists 
of \emph{non-commutative polynomials} $p$  of the form
$$p(z)=\sum_{w\in\Lambda}p_wz^w,$$
where $\Lambda\subset\mathcal{F}_N$ is a finite set. We will consider also the space $\mathcal{L(U,Y)}\left\langle \left\langle z_1,\ldots,z_N\right\rangle\right\rangle$ of formal power series with coefficients in $\mathcal{L(U,Y)}$ and the space $\mathcal{L(U,Y)} \left\langle z_1,\ldots,z_N\right\rangle$ of non-commutative polynomials with coefficients in $\mathcal{L(U,Y)}$.

Let us introduce now the \emph{non-commutative Herglotz--Agler class} $\mathcal{HA}_N^{\rm nc}(\mathcal{Y})$ of formal power series $f\in\mathcal{L(Y)}\left\langle \left\langle z_1,\ldots,z_N\right\rangle\right\rangle$ such that 
the series
$$f(\mathbf{C}):=\sum_{w\in\mathcal{F}_N}f_w\otimes\mathbf{C}^w$$
converges in the operator norm and
$\Re f(\mathbf{C})\geq 0$ for every $\mathbf{C}\in\mathcal{D}^N\cap\mathcal{L(E)}^N$, with a Hilbert space $\mathcal{E}$. Here for $w=g_{i_1}\cdots g_{i_m}\in\mathcal{F}_N$ we set $\mathbf{C}^w:=C_{i_1}\cdots C_{i_m}$, and $\mathbf{C}^\emptyset =I_\mathcal{E}$. The subclass $\mathcal{HA}_N^{{\rm nc},I}(\mathcal{Y})\subset\mathcal{HA}_N^{\rm nc}(\mathcal{Y})$ consists of formal power series $f$ such that $f_\emptyset =I_\mathcal{Y}$.
\begin{thm}\label{thm:ha-nc}
A formal power series $f\in\mathcal{L(Y)}\left\langle \left\langle z_1,\ldots,z_N\right\rangle\right\rangle$ belongs to the class $\mathcal{HA}^{nc,I}_N(\mathcal{Y})$ if and only if there exist a Hilbert space $\mathcal{H}$, an $N$-tuple $\mathbf{G}\in\mathcal{G}_N\cap\mathcal{L(H)}^N$, and an isometry $V\in\mathcal{L(Y,H)}$ such that
\begin{equation}\label{ha-nc}
f(z)=V^*(I_\mathcal{H}+z\mathbf{G})(I_\mathcal{H}-z\mathbf{G})^{-1}V,
\end{equation}
where $z\mathbf{G}:=\sum_{k=1}^Nz_kG_k$ and thus
\begin{equation}\label{res}
(I_\mathcal{H}-z\mathbf{G})^{-1}=\sum_{j=0}^\infty\left(\sum_{k=1}^Nz_kG_k\right)^j=\sum_{w\in\mathcal{F}_N}\mathbf{G}^wz^w.
\end{equation}
\end{thm}
\begin{proof}
If \eqref{ha-nc} holds then $f\in\mathcal{HA}^{nc,I}_N(\mathcal{Y})$. Indeed, for any $\mathbf{C}\in\mathcal{D}^N\cap\mathcal{L(E)}^N$, with a Hilbert space $\mathcal{E}$, by Corollary~\ref{cor:gn} the operator  $\mathbf{G\otimes C}$ is a strict contraction. Then the series in \eqref{res} evaluated at $\mathbf{C}$ converges in the operator norm, thus the operator $(I_\mathcal{H\otimes E}+\mathbf{G\otimes C})(I_\mathcal{H\otimes E}-\mathbf{G\otimes C})^{-1}$ is well defined and, as the Cayley transform of a strict contraction, has positive semidefinite real part, and so is $f(\mathbf{C})$. Clearly, since $V$ is an isometry, $f_\emptyset=I_\mathcal{Y}$.

Let us prove the converse.
The formal power series 
\begin{equation}\label{cayley}
F(z):=(f(z)-I_\mathcal{Y})(f(z)+I_\mathcal{Y})^{-1}
\end{equation}
is well defined since $f_\emptyset+I_\mathcal{Y}=2I_\mathcal{Y}$ is invertible and so is $f(z)+I_\mathcal{Y}$. Moreover, $F$ belongs to the \emph{non-commutative Schur--Agler class $\mathcal{SA}_N^{\rm nc}(\mathcal{Y})$}, i.e., the formal power series  $F$ evaluated at any $\mathbf{C}\in\mathcal{D}^N$ is well defined and $\| F(\mathbf{C})\|\leq 1$. By \cite{BGM}, there exists a Hilbert space $\mathcal{H}$, a resolution of identity $\mathbf{P}=(P_1,\ldots,P_N)\in\mathcal{L(H)}^N$, and a unitary operator {\small $U=\left[\begin{array}{cc} A & B\\ C & D\end{array}\right]\in\mathcal{L(H\oplus Y)}$} such that
\begin{equation}\label{sa-nc}
F(z)=D+C(I_\mathcal{Y}-(z\mathbf{P})A)^{-1}(z\mathbf{P})B
\end{equation}
where we set $z\mathbf{P}:=\sum_{k=1}^Nz_kP_k$. To get the representation \eqref{ha-nc} for $f$, we will apply a trick which is well known in one-variable system theory (see, e.g., \cite{Ar}). Consider a \emph{non-commutative linear system}
$\Sigma =(N; U; \mathbf{P}; \mathcal{H,Y})$, i.e., a system of equations
\begin{equation}\label{scat}
\left\{ \begin{array}{lll}
x(z) &=& (z\mathbf{P})Ax(z)+(z\mathbf{P})Bu(z),\\
y(z) &=& Cx(z)+Du(z),
\end{array}\right.
\end{equation}
where $x(z)\in\mathcal{L(H)}\left\langle \left\langle z_1,\ldots,z_N\right\rangle\right\rangle,\ u(z),y(z)\in\mathcal{L(Y)}\left\langle \left\langle z_1,\ldots,z_N\right\rangle\right\rangle $ (the corresponding system of equations for coefficients of these formal power series is one of the systems with evolution on the free semigroup considered in \cite{BGM}). The system \eqref{scat} is equivalent to the system
\begin{equation}\label{scat'}
\left\{ \begin{array}{lll}
x(z) &=& (I_\mathcal{H}-(z\mathbf{P})A)^{-1}(z\mathbf{P})Bu(z),\\
y(z) &=& F(z)u(z).
\end{array}\right.
\end{equation}
The second equation means that $F$ is the \emph{transfer function} of the system $\Sigma$, or $\Sigma$ is a \emph{realization} of the formal power series $F$. Let us find a realization $\widetilde{\Sigma} =(N; \widetilde{U}; \mathbf{P}; \mathcal{H,Y})$ of $f$. To this end (now the trick appears!) we apply the so-called \emph{diagonal transform}:
\begin{equation}\label{diag}
u(z)=\widetilde{y}(z)+\widetilde{u}(z),\qquad y(z)=\widetilde{y}(z)-\widetilde{u}(z).
\end{equation}
Then we get $\widetilde{y}(z)=f(z)\widetilde{u}(z)$, i.e., an analogue of the second equation in \eqref{scat'}. Suppose that the operator $I_\mathcal{Y}-D$ is invertible. Then an easy calculation gives the desired system realization $\widetilde{\Sigma}$:
\begin{equation}\label{realiz}
\left\{ \begin{array}{lll}
x(z) &=& (z\mathbf{P})\widetilde{A}x(z)+(z\mathbf{P})\widetilde{B}u(z),\\
y(z) &=& \widetilde{C}x(z)+\widetilde{D}u(z),
\end{array}\right.
\end{equation}
where
\begin{eqnarray*}
\widetilde{A}=A+B(I_\mathcal{Y}-D)^{-1}C, & \widetilde{B}=2B(I_\mathcal{Y}-D)^{-1}\\
\widetilde{C}=(I_\mathcal{Y}-D)^{-1}C,   & \widetilde{D}=(I_\mathcal{Y}-D)^{-1}(I_\mathcal{Y}+D).
\end{eqnarray*}
In our case $\widetilde{D}=f_\emptyset=I_\mathcal{Y}$ and $D=F_\emptyset=0$. Then
$$
\widetilde{A}=A+BC, \quad \widetilde{B}=2B,\quad
\widetilde{C}=C.
$$
Moreover, since $U$ is a unitary operator, in this case $B$ is an isometry, $C$ is a coisometry, $A+BC$ is unitary, and $A^*B=0,\ AC^*=0$. Thus, we may write
\begin{eqnarray*}
f(z) &=& I_\mathcal{Y}+2C(I_\mathcal{H}-(z\mathbf{P})(A+BC))^{-1}(z\mathbf{P})B\\
     &=& I_\mathcal{Y}+2C(I_\mathcal{H}-(z\mathbf{P})(A+BC))^{-1}(z\mathbf{P})(A+BC)(A+BC)^*B\\
     &=& I_\mathcal{Y}+2C(I_\mathcal{H}-(z\mathbf{P})(A+BC))^{-1}(z\mathbf{P})(A+BC)C^*\\
     &=& C(I_\mathcal{H}+(z\mathbf{P})(A+BC))(I_\mathcal{H}-(z\mathbf{P})(A+BC))^{-1}C^*\\
     &=& V^*(I_\mathcal{H}+z\mathbf{G})(I_\mathcal{H}-z\mathbf{G})^{-1}V,
\end{eqnarray*}
where $G_k=P_k(A+BC),\ k=1,\ldots,N$ and $V=C^*$. By Proposition~\ref{prop:gn}, $\mathbf{G}\in\mathcal{G}_N$. Since $C$ is a coisometry, $V=C^*$ is an isometry. Thus we have obtained a representation \eqref{ha-nc} of $f$.
\end{proof}
\begin{cor}\label{cor:coef}
A formal power series $f\in\mathcal{L(Y)}\left\langle \left\langle z_1,\ldots,z_N\right\rangle\right\rangle$ belongs to the class $\mathcal{HA}_N^{nc}(\mathcal{Y})$ and satisfies $f_\emptyset=\frac{I_\mathcal{Y}}{2}$ if and only if there exist a Hilbert space $\mathcal{H}$, an $N$-tuple of operators $\mathbf{G}\in\mathcal{G}_N\cap\mathcal{L(H)}^N$, and an isometry $V\in\mathcal{L(Y,H)}$ such that the sequence 
$$
a_\emptyset:=I_\mathcal{Y},\ a_w:=f_w \qquad {\rm for}\quad w\in\mathcal{F}_N\setminus\{\emptyset\}, 
$$
satisfies
$$
a_w=V^*\mathbf{G}^wV,\qquad w\in\mathcal{F}_N. 
$$
\end{cor}
\begin{proof}
The statement follows from the representation \eqref{ha-nc} for $\widetilde{f}=2f$:
\begin{eqnarray*}
\widetilde{f}(z) &=& V^*(2(I_\mathcal{H}-z\mathbf{G})^{-1}-I_\mathcal{H})V=V^*\left(2\sum_{j=0}^\infty (z\mathbf{G})^j-I_\mathcal{H}\right)V \\
              &=& I_\mathcal{H}+2V^*\left(\sum_{w\in\mathcal{F}_N\setminus\{\emptyset\}}\mathbf{G}^wz^w\right)V.
\end{eqnarray*}
\end{proof}
\begin{rem}\label{rem:mpd}
In \cite{AK2} it has been shown that a formal power series $F$ belongs to the class $\mathcal{SA}_N^{\rm nc}(\mathcal{Y})$ (or more generally, to the class $\mathcal{SA}_N^{\rm nc}(\mathcal{U,Y})$) if and only if the series for $F(\mathbf{C})$ converges to a contractive operator for every $\mathbf{C}\in\mathcal{D}^N_{\rm matr}$, i.e., it is enough to test values of $F$ on $N$-tuples of strictly contractive matrices of same size $n\times n,\ n=1,2,\ldots$. The analogous statement for the class $\mathcal{HA}_N^{\rm nc}(\mathcal{Y})$ is true, too: a formal power series $f$ belongs to the class $\mathcal{HA}_N^{\rm nc}(\mathcal{Y})$ if and only if the series for $f(\mathbf{C})$ converges and $\Re f(\mathbf{C})$ is positive semidefinite for every $\mathbf{C}\in\mathcal{D}^N_{\rm matr}$. Indeed, this follows from the fact that the Cayley transform $f\mapsto F$ defined by \eqref{cayley} is an injection from $\mathcal{HA}_N^{\rm nc}(\mathcal{Y})$ into $\mathcal{SA}_N^{\rm nc}(\mathcal{Y})$. 
\end{rem}

\section{The Carath\'{e}odory interpolation problem}\label{sec:prob}
In this section we will consider $\mathcal{F}_N$ as a sub-semigroup of the \emph{free semigroup with involution} $\hat{\mathcal{F}}_{2N}$. The latter is the free semigroup $\mathcal{F}_{2N}$ with the generators $g_1,\ldots,g_N,g_{N+1},\ldots,g_{2N}$ and the neutral element $\emptyset$, endowed with the involution ``$*$" defined as follows: $g_k^*:=g_{k+N}$ for $k=1,\ldots,N$, $g_k^*:=g_{k-N}$ for $k=N+1,\ldots,2N$, $\emptyset^*:=\emptyset$, and $(g_{i_1}\cdots g_{i_l})^*:=g_{i_l}^*\cdots g_{i_1}^*$ for every $l\in\mathbb{N}$ and $i_j\in\{1,\ldots,2N\},\ j=1,\ldots,l$. For a set $\Omega\subset\hat{\mathcal{F}}_{2N}$ we define the set $\Omega^*:=\{ w\in\hat{\mathcal{F}}_{2N}:\ w^*\in\Omega\}.$ Let us introduce also the unital $*$-algebra $\mathcal{A}_N(\mathcal{Y})$ as the algebra $\mathcal{L(Y)}\left\langle z_1,\ldots,z_N,z_{N+1},\ldots z_{2N}\right\rangle$  endowed with the involution ``$*$" defined as follows: 1) $z_k^*:=z_{k+N}$ for $k=1,\ldots,N$, $z_k^*:=z_{k-N}$ for $k=N+1,\ldots,2N$, $(z_{i_1}\cdots z_{i_l})^*:=z_{i_l}^*\cdots z_{i_1}^*$ for every $l\in\mathbb{N}$ and $i_j\in\{1,\ldots,2N\},\ j=1,\ldots,l$, thus for $\hat{z}:=(z_1,\ldots,z_N,z_{N+1},\ldots z_{2N})=(z_1,\ldots,z_N,z_1^*,\ldots,z_N^*)$ and $w\in\hat{\mathcal{F}}_{2N}$ one has $(\hat{z}^w)^*=\hat{z}^{w^*}$; 2) for arbitrary finite set $\Omega\subset\hat{\mathcal{F}}_{2N}$ and a polynomial $p(\hat{z})=\sum_{w\in\Omega}p_w\hat{z}^w$, $$p(\hat{z})^*=\left(\sum_{w\in\Omega}p_w\hat{z}^w\right)^*:=\sum_{w\in\Omega}p_w^*\hat{z}^{w^*}=\sum_{w\in\Omega^*}p_{w^*}^*\hat{z}^w$$
(here $p_w^*$ is the adjoint operator to $p_w$ in $\mathcal{L(Y)}$). 

A finite set $\Lambda\subset\mathcal{F}_N$ will be called \emph{admissible} if $g_kw\in\mathcal{F}_N\setminus\Lambda$ and $wg_k\in\mathcal{F}_N\setminus\Lambda$ for every $w\in\mathcal{F}_N\setminus\Lambda$ and $k=1,\ldots,N$. Clearly, if the set $\Lambda$ is admissible and non-empty then $\emptyset\in\Lambda$, and if $\Lambda$ is admissible, non-empty and $\Lambda\neq\{\emptyset\}$ then there is a $k\in\{ 1,\ldots,N\}$ such that $g_k\in\Lambda$. For example, the set $\Lambda_m:=\{ w\in\mathcal{F}_N:\ |w|\leq m\}$ is admissible.

Let us pose now the \emph{Carath\'{e}odory interpolation problem in the class $\mathcal{HA}_N^{\rm nc}(\mathcal{Y})$}.
\begin{prob}\label{prob:c-gen}
Let $\Lambda\subset\mathcal{F}_N$ be an admissible set. Given a collection of operators $\{ c_w\}_{w\in\Lambda}\in\mathcal{L(Y)}$, with $c_\emptyset\geq 0$, find 
$f\in\mathcal{HA}_N^{\rm nc}(\mathcal{Y})$ such that
$$
f_\emptyset=\frac{c_\emptyset}{2},\quad f_w=c_w\quad {\rm for}\ w\in\Lambda\setminus\{\emptyset\}.
$$
\end{prob}
We will start with the special case of this problem where $c_\emptyset=I_\mathcal{Y}$.
\begin{prob}\label{prob:c}
Let $\Lambda\subset\mathcal{F}_N$ be an admissible set. Given a collection of operators $\{ c_w\}_{w\in\Lambda}\in\mathcal{L(Y)}$, with $c_\emptyset=I_\mathcal{Y}$, find 
$f\in\mathcal{HA}_N^{\rm nc}(\mathcal{Y})$ such that
$$f_\emptyset=\frac{c_\emptyset}{2}=\frac{I_\mathcal{Y}}{2},\quad f_w=c_w\quad {\rm for}\ w\in\Lambda\setminus\{\emptyset\}.
$$
\end{prob}
 
From Corollary~\ref{cor:coef} we obtain the following result.
\begin{thm}\label{thm:g}
Problem~\ref{prob:c} has a solution if and only if there exist a Hilbert space $\mathcal{H}$, an $N$-tuple of operators $\mathbf{G}\in\mathcal{G}_N\cap\mathcal{L(H)}^N$, and an isometry $V\in\mathcal{L(Y,H)}$ such that
\begin{equation}\label{mom-nc}
c_w=V^*\mathbf{G}^wV,\qquad w\in\Lambda.
\end{equation}
\end{thm}
\begin{rem}\label{rem:arb}
Theorem~\ref{thm:g} holds true for the case where Problem~\ref{prob:c} is formulated for an arbitrary set $\Lambda\subset\mathcal{F}_N$, not necessarily finite and admissible.
\end{rem}
Let us note that \eqref{mom-nc} is a non-commutative multivariable counterpart of \eqref{mom-o}. Theorem~\ref{thm:g} gives a criterion on solvability of Problem~\ref{prob:c} in the ``existence terms". We are going to obtain also another criterion, in terms of positivity of certain non-commutative polynomial whose coefficients are determined by the problem data. 

Let $\mathcal{E}$ be a Hilbert space, and $\mathbf{T}=(T_1,\ldots,T_N)\in\mathcal{L(E)}^N$. Then set $$\hat{\mathbf{T}}:=(T_1,\ldots,T_N,T_1^*,\ldots,T_N^*)\in\mathcal{L(E)}^{2N}.$$ For a finite set $\Omega\subset\hat{\mathcal{F}}_{2N}$ and a polynomial $p(\hat{z})=\sum_{w\in\Omega}p_w\hat{z}^w\in\mathcal{A}_N(\mathcal{Y})$ define $$p(\mathbf{T})=p(\hat{\mathbf{T}}):=\sum_{w\in\Omega}p_w\otimes\hat{\mathbf{T}}^w\in\mathcal{L(Y\otimes E)}.$$
In particular, if $\Omega=\Lambda\cup\Lambda^*$ where $\Lambda\subset\mathcal{F}_N$ is a finite set, and $$p(\hat{z})=\sum_{w\in\Lambda\cup\Lambda^*}p_w\hat{z}^w=p_\emptyset+\sum_{w\in\Lambda\setminus\{\emptyset\}}p_wz^w+\sum_{w\in\Lambda^*\setminus\{\emptyset\}}p_wz^{*w},$$ we have 
$$p(\mathbf{T})=\sum_{w\in\Lambda\cup\Lambda^*}p_w\otimes\hat{\mathbf{T}}^w=p_\emptyset\otimes I_\mathcal{E}+\sum_{w\in\Lambda\setminus\{\emptyset\}}p_w\otimes \mathbf{T}^w+\sum_{w\in\Lambda^*\setminus\{\emptyset\}}p_w\otimes\mathbf{T}^{*w},$$
where $z^*:=(z_1^*,\ldots,z_N^*)$ and $\mathbf{T}^*:=(T_1^*,\ldots,T_N^*)$, and one identifies $z^w=\hat{z}^w,\ \mathbf{T}^w=\hat{\mathbf{T}}^w$ for $w\in\mathcal{F}_N\subset\hat{\mathcal{F}}_{2N}$, and $z^{*w}=\hat{z}^{*w},\ \mathbf{T}^{*w}=\hat{\mathbf{T}}^{*w}$ for $w\in\mathcal{F}_N^*\subset\hat{\mathcal{F}}_{2N}$.
Thus, the evaluation of polynomials from $\mathcal{A}_N(\mathcal{Y})$ on $N$-tuples of operators is well defined.  
\begin{lem}\label{lem:main}
Let $\emptyset\in\Lambda\subset\mathcal{F}_N$ be a finite set.
\begin{description}
	\item[I] A polynomial $$p(\hat{z})=I_\mathcal{Y}+\sum_{w\in\Lambda\setminus\{\emptyset\}}p_wz^w+\sum_{w\in\Lambda^*\setminus\{\emptyset\}}p_wz^{*w}\in\mathcal{A}_N(\mathcal{Y})$$ is positive semidefinite on $\mathcal{U}^N$ if and only if there exist a Hilbert space $\mathcal{H}$, an $N$-tuple of operators $\mathbf{G}\in\mathcal{G}_N\cap\mathcal{L(H)}^N$, and an isometry $V\in\mathcal{L(Y,H)}$ such that 
	\begin{eqnarray}
	p_w=p_{w^*}^* &=& V^*\mathbf{G}^wV,\qquad w\in\Lambda,\label{I1}\\
	0 &=& V^*\mathbf{G}^wV,\qquad w\in \mathcal{F}_N\setminus \Lambda.\label{I2}
\end{eqnarray}	
\item[II] A polynomial $$p(\hat{z})=I_\mathcal{Y}+\sum_{w\in\Lambda\setminus\{\emptyset\}}p_wz^w+\sum_{w\in\Lambda^*\setminus\{\emptyset\}}p_wz^{*w}\in\mathcal{A}_N(\mathcal{Y})$$ is positive semidefinite on $\mathcal{G}^N$ if and only if there exist a Hilbert space $\mathcal{K}$, an $N$-tuple of operators $\mathbf{U}\in\mathcal{U}^N\cap\mathcal{L(K)}^N$, and an isometry $W\in\mathcal{L(Y,K)}$ such that 
	\begin{eqnarray}
	p_w=p_{w^*}^* &=& W^*\mathbf{U}^wW,\qquad w\in\Lambda,\label{II1}\\
	0 &=& W^*\mathbf{U}^wW,\qquad w\in \mathcal{F}_N\setminus \Lambda.\label{II2}
\end{eqnarray}
\end{description}
\end{lem}
\begin{proof}
I. If the polynomial $p$ is positive semidefinite on $\mathcal{U}^N$ then $p_w=p_{w^*}^*$ for $w\in\Lambda$. This can be seen from a McCullough factorization \cite{McC3}: $p(\hat{z})=h(z)^*h(z)$, where $h(z)=\sum_{w\in\mathcal{F}_N:\,|w|\leq m}h_wz^w\in\mathcal{L(Y,V)}\left\langle z_1,\ldots,z_N\right\rangle$, with an auxiliary Hilbert space $\mathcal{V}$. Set $$f(z):=\frac{I_\mathcal{Y}}{2}+\sum_{w\in\Lambda\setminus\{\emptyset\}}p_wz^w.$$
Then $p(\mathbf{C})=2\Re f(\mathbf{C})\geq 0$ for every $\mathbf{C}\in\mathcal{D}^N\cap\mathcal{L(E)}^N$, with a Hilbert space $\mathcal{E}$. Indeed, since $\mathbf{C}$ has a unitary dilation (see \cite{SzNF}) $\mathbf{U}\in\mathcal{U}^N\cap\mathcal{L(K)}^N$, i.e., $\mathcal{K}\supset\mathcal{E}$ and
$$\mathbf{C}^w=P_\mathcal{E}\mathbf{U}^w\big|_{\mathcal{E}},\qquad w\in\mathcal{F}_N,$$
and $p(\mathbf{U})=2\Re f(\mathbf{U})\geq 0$, we get
$$p(\mathbf{C})=(I_\mathcal{Y}\otimes P_\mathcal{E})p(\mathbf{U})\big|_\mathcal{Y\otimes E}\geq 0.$$ 
Thus, $f\in\mathcal{HA}^{\rm nc}_N(\mathcal{Y})$ and $f_\emptyset = \frac{I_\mathcal{Y}}{2}$. By Corollary~\ref{cor:coef}, there exists a representation \eqref{I1}--\eqref{I2}.

Conversely, if $p$ has a representation \eqref{I1}--\eqref{I2} then $p(\hat{z})= f(z)+f(z)^*$, where $$f(z)=\frac{I_\mathcal{Y}}{2}+\sum_{w\in\Lambda\setminus\{\emptyset\}}p_wz^w=\frac{I_\mathcal{Y}}{2}+\sum_{w\in\Lambda\setminus\{\emptyset\}}V^*\mathbf{G}^wVz^w.$$
By Corollary~\ref{cor:coef}, $f\in\mathcal{HA}^{\rm nc}_N(\mathcal{Y})$. Hence, $$p(\mathbf{U})=2\Re f(\mathbf{U})=\lim_{r\uparrow 1}2\Re f(r\mathbf{U})\geq 0,\qquad \mathbf{U}\in\mathcal{U}^N.$$

II. Let $p$ be positive semidefinite on $\mathcal{G}_N$. Let $\mathcal{A}_{\mathcal{U}^N}$ be the $C^*$-algebra obtained as  the norm completion of the quotient of unital $*$-algebra $\mathcal{A}_N=\mathcal{A}_N(\mathbb{C})$ with the seminorm
$$\| q\| :=\sup_{\mathbf{U}\in\mathcal{U}^N}\| q(\mathbf{U})\|=\sup_{\mathbf{U}\in\mathcal{U}^N}\| q(U_1,\ldots,U_N,U_1^*,\ldots,U_N^*)\|,$$
by the two-sided ideal of elements of zero seminorm. 

Let us show that the restriction of the quotient map above to the subspace $\mathcal{B}_N\subset\mathcal{A}_N$ of polynomials of the form 
\begin{equation}\label{q}
q(\hat{z})=q_\emptyset +\sum_{w\in\mathcal{F}_N:\,0<|w|\leq m}q_wz^w+\sum_{w\in\mathcal{F}_N^*:\,0<|w|\leq m}q_wz^{*w}
\end{equation}
 is injective, i.e., that if $q\in\mathcal{B}_N$ is non-zero then the corresponding coset $[q]\in\mathcal{A}_{\mathcal{U}^N}$ is non-zero. Indeed, if $[q]=[0]$ then $q(\mathbf{U})=0$ for every $\mathbf{U}\in\mathcal{U}^N$. In particular, $q$ is positive semidefinite on $\mathcal{U}^N$. By \cite{McC3}, there exists a polynomial $h(z)=\sum_{w\in\mathcal{F}_N:\,|w|\leq m}h_wz^w\in\mathcal{L}(\mathbb{C},\mathbb{C}^r)\left\langle z_1,\ldots,z_N\right\rangle$, with some $r\in\mathbb{N}$, such that $q(\hat{z})=h(z)^*h(z)$. Then 
$h$ vanishes on $\mathcal{U}^N$. In particular, for every $n\in\mathbb{N}$ the polynomial $h$ vanishes on $\mathcal{U}^N\cap (\mathbb{C}^{n\times n})^N$. The latter set is the uniqueness set for functions of matrix entries, which are holomorphic on a domain containing this set (see, e.g., \cite{Sh}). Then $h$ vanishes on the all of $(\mathbb{C}^{n\times n})^N$, for each $n\in\mathbb{N}$. By the Amitsur--Levitzki theorem (see \cite[pp. 22-23]{Ro}) such a polynomial should be zero, i.e., $h_w=0$ for all $w\in\mathcal{F}_N:\,|w|\leq m$. Then $q_w=0$ for all $w\in\mathcal{F}_N\cup\mathcal{F}_N^*:\,|w|\leq m$. 

Denote by $\mathcal{B}_{\mathcal{U}^N}$ the image of the subspace $\mathcal{B}_N$ under the quotient map above. This subspace of the $C^*$-algebra $\mathcal{A}_{\mathcal{U}^N}$ is \emph{selfadjoint}, i.e., $$\mathcal{B}_{\mathcal{U}^N}^*:=\{ [q]^*=[q^*]:\ [q]\in\mathcal{B}_{\mathcal{U}^N}\}=\mathcal{B}_{\mathcal{U}^N}.$$ Define the linear map $\varphi:\ \mathcal{B}_{\mathcal{U}^N}\to\mathcal{L(Y)}$ by  $\varphi([z^w])=p_w$ for $w\in\Lambda$, $\varphi([z^{*w}])=p_w$ for $w\in\Lambda^*$, and $\varphi([z^w])=\varphi([z^{*w^*}])=0$ for $w\in\mathcal{F}_N\setminus\Lambda$. By the result of the previous paragraph together with the Amitsur--Levitzki theorem mentioned there, this linear map is correctly defined. Let us show that $\varphi$ is \emph{completely positive}, i.e., that for every $n\in\mathbb{N}$ the map $$\varphi_n:= {\rm id}_n\otimes\varphi:\ \mathbb{C}^{n\times n} \otimes\mathcal{B}_{\mathcal{U}^N}\to\mathbb{C}^{n\times n}\otimes\mathcal{L(Y)}$$ (here ${\rm id}_n$ is the identity map from the $C^*$-algebra $\mathbb{C}^{n\times n}$ onto itself) is \emph{positive}. The latter means, in turn, that $\varphi_n$ maps positive elements (in the sense  of the $C^*$-algebra $\mathbb{C}^{n\times n}\otimes\mathcal{A}_{\mathcal{U}^N}$) from $\mathbb{C}^{n\times n}\otimes\mathcal{B}_{\mathcal{U}^N}$ into positive elements in the $C^*$-algebra $\mathbb{C}^{n\times n}\otimes\mathcal{L(Y)}$. Let $[q]\in\mathbb{C}^{n\times n}\otimes\mathcal{B}_{\mathcal{U}^N}$ be a positive element of the $C^*$-algebra $\mathbb{C}^{n\times n}\otimes\mathcal{A}_{\mathcal{U}^N}$, i.e.,
$[q]=[h]^*[h]$ with some $[h]\in\mathbb{C}^{n\times n}\otimes\mathcal{A}_{\mathcal{U}^N}$. One can think of $[q]$ as of the $n\times n$ matrix $([q]_{ij})_{i,j=1,\ldots,n}$ whose entries $[q]_{ij}=[q_{ij}]\in\mathcal{B}_{\mathcal{U}^N}$ and $q_{ij}\in\mathcal{B}_N$, and thus $q\in\mathbb{C}^{n\times n}\otimes\mathcal{B}_N$ is a polynomial of the form \eqref{q} with the coefficients from $\mathbb{C}^{n\times n}$. Let us observe that by virtue of the definition of the $C^*$-algebra $\mathcal{A}_{\mathcal{U}^N}$, for an arbitrary $[x]\in\mathcal{A}_{\mathcal{U}^N}$ its values on $\mathcal{U}^N$ are well defined. In particular, if $x\in\mathcal{B}_N$ then $[x](\mathbf{U})=x(\mathbf{U})$ for any $\mathbf{U}\in\mathcal{U}^N$. Therefore, for an arbitrary $[x]=([x]_{ij})_{i,j=1,\ldots,n}=([x_{ij}])_{i,j=1,\ldots,n}\in\mathbb{C}^{n\times n}\otimes\mathcal{A}_{\mathcal{U}^N}$ one defines correctly $[x](\mathbf{U}):=([x_{ij}](\mathbf{U}))_{i,j=1,\ldots,n},\ \mathbf{U}\in\mathcal{U}^N$. In particular, if $x=(x_{ij})_{i,j=1,\ldots,n}\in\mathbb{C}^{n\times n}\otimes\mathcal{B}_N$ then $[x](\mathbf{U})=x(\mathbf{U})=(x_{ij}(\mathbf{U}))_{i,j=1,\ldots,n}$ for any $\mathbf{U}\in\mathcal{U}^N$. Since $q(\mathbf{U})=[q](\mathbf{U})=[h](\mathbf{U})^*[h](\mathbf{U})$ is positive semidefinite for every  $\mathbf{U}\in\mathcal{U}^N$, it follows that the polynomial $q$ (with the coefficients from $\mathbb{C}^{n\times n}$) is positive semidefinite on $\mathcal{U}^N$. From the McCullough factorization theorem \cite{McC3} we deduce that $q_\emptyset\geq 0$. If $q_\emptyset=I_n$, then by Part I of this Lemma, there exist a Hilbert space $\mathcal{H}$, an $N$-tuple of operators $\mathbf{G}\in\mathcal{G}_N\cap\mathcal{L}(\mathcal{H})^N$, and an isometry $V\in\mathcal{L}(\mathbb{C}^n,\mathcal{H})$ such that 
\begin{eqnarray*}
q_w=q_{w^*}^* &=& V^*\mathbf{G}^wV,\qquad w\in\mathcal{F}_N:\,|w|\leq m,\\
0   &=& V^*\mathbf{G}^wV,\qquad w\in\mathcal{F}_N:\,|w|>m.
\end{eqnarray*}
Then we have
\allowdisplaybreaks
\begin{eqnarray*}
\varphi_n([q]) &=& ({\rm id}_n\otimes \varphi)\left(I_n\otimes [1]+\sum_{w\in\mathcal{F}_N:\,0<|w|\leq m}q_w\otimes [z^w]\right.\\
&+& \left.\sum_{w\in\mathcal{F}_N^*:\,0<|w|\leq m}q_w\otimes [z^{*w}] \right)\\ 
&=& I_{\mathbb{C}^n\otimes \mathcal{Y}}+\sum_{w\in\Lambda\setminus\{\emptyset\} }V^*\mathbf{G}^wV\otimes p_w+ \sum_{w\in\Lambda^*\setminus\{\emptyset\} }V^*\mathbf{G}^{*w}V\otimes p_w\\
&=& (V^*\otimes I_\mathcal{Y})p^{\rm l}(\mathbf{G})(V\otimes I_\mathcal{Y})\geq 0
\end{eqnarray*}
(positivity in the $C^*$-algebra $\mathbb{C}^{n\times n}\otimes\mathcal{L(Y)}\cong\mathcal{L}(\mathbb{C}^n\otimes\mathcal{Y})$ is operator positive semidefiniteness). 
In the case where $q_\emptyset >0$ we  define $\tilde{q}(\hat{z}):=q_\emptyset^{-1/2}q(\hat{z})q_\emptyset^{-1/2}$. Since $\varphi_n([\tilde{q}])\geq 0$, we get $$\varphi_n([q])=(q_\emptyset^{1/2}\otimes I_\mathcal{Y})\varphi_n([\tilde{q}])(q_\emptyset^{1/2}\otimes I_\mathcal{Y})\geq 0.$$
In the case where the matrix $q_\emptyset$ is degenerate we set $q_\epsilon(\hat{z}) :=\epsilon I_n+q(\hat{z})$ for  $\epsilon >0$. Then $q_\epsilon$ is positive definite on $\mathcal{U}^N$ and ${(q_\epsilon)}_\emptyset=\epsilon I_n+q_\emptyset >0$. Since $\varphi_n([q_\epsilon])\geq 0$, we get
 $$
\varphi_n([q]) = \lim_{\epsilon\downarrow 0}\varphi_n([q_\epsilon ])\geq 0.
$$
Finally, we have obtained that $\varphi: \mathcal{B}_{\mathcal{U}^N}\rightarrow \mathcal{L(Y)}$ is completely positive.

Since we have $\varphi ([1])=I_\mathcal{Y}$, by the Arveson extension theorem \cite{Arv1} there exists a completely positive map $\widetilde{\varphi}:\mathcal{A}_{\mathcal{U}^N}\rightarrow \mathcal{L(Y)}$ which extends $\varphi$. By the Stinespring theorem \cite{St}, there exists a $*$-representation $\pi$ of $\mathcal{A}_{\mathcal{U}^N}$ in some Hilbert space $\mathcal{K}$ and an isometry $W\in\mathcal{L(Y,K)}$ such that
$$\widetilde{\varphi}(a)=W^*\pi(a)W,\qquad a\in\mathcal{A}_{\mathcal{U}^N}.$$
In particular, we get
\begin{eqnarray}
p_w &=& \varphi([z^w])=\widetilde{\varphi}([z^w])=W^*\pi([z^w])W=W^*\mathbf{U}^wW,\quad w\in\Lambda,\label{II1'}\\
p_w &=& \varphi([z^{*w}])=\widetilde{\varphi}([z^{*w}])=W^*\pi([z^{*w}])W=W^*\mathbf{U}^{*w}W,\quad w\in\Lambda^*,\label{II2'}\\
0   &=& \varphi([z^w])=\widetilde{\varphi}([z^w])=W^*\pi([z^w])W=W^*\mathbf{U}^wW,\quad w\in \mathcal{F}_N\setminus \Lambda,
\label{II3'}
\end{eqnarray}
where we set $\mathbf{U}:=(\pi([z_1]),\ldots,\pi([z_N]))$. We have $\mathbf{U}\in\mathcal{U}^N$. Indeed, since $$\| [1-z_k^*z_k]\|_{\mathcal{A}_{\mathcal{U}^N}}=\| [1-z_kz_k^*]\|_{\mathcal{A}_{\mathcal{U}^N}}=0,$$ we get $[z_k^*z_k]=[z_kz_k^*]=[1]$. Hence
\begin{eqnarray*} 
U_k^*U_k &=& \pi([z_k])^*\pi([z_k])=\pi([z_k^*z_k])=\pi([1])=I_\mathcal{K},\quad k=1,\ldots,N,\\ U_kU_k^* &=& \pi([z_k])\pi([z_k])^*=\pi([z_kz_k^*])=\pi([1])=I_\mathcal{K},\quad k=1,\ldots,N.
\end{eqnarray*}
Clearly, \eqref{II1'} and \eqref{II2'} imply that $p_w=p_{w^*}^*$ for $w\in\Lambda$.
Thus, representation \eqref{II1'}--\eqref{II3'} for the coefficients of $p$ is a desired representation \eqref{II1}--\eqref{II2}.

Conversely, if the coefficients of $p$ have a representation \eqref{II1}--\eqref{II2} then for any $\mathbf{G}\in\mathcal{G}_N\cap\mathcal{L}(\mathcal{H})^N$ one has
\begin{eqnarray*}
p(\mathbf{G}) &=& I_\mathcal{Y\otimes H}+ \sum_{w\in\Lambda\setminus\{\emptyset\} }p_w\otimes\mathbf{G}^w+\sum_{w\in\Lambda^*\setminus\{\emptyset\}}p_w\otimes\mathbf{G}^{*w}\\
&=& I_\mathcal{Y\otimes H}+ \sum_{w\in\Lambda\setminus\{\emptyset\}}W^*\mathbf{U}^wW\otimes\mathbf{G}^w+\sum_{w\in\Lambda^*\setminus\{\emptyset\}}W^*\mathbf{U}^{*w}W\otimes\mathbf{G}^{*w}\\
&=& I_\mathcal{Y\otimes H}+ \sum_{w\in\mathcal{F}_N\setminus\{\emptyset\}}W^*(\mathbf{U})^wW\otimes\mathbf{G}^w+
\sum_{w\in\mathcal{F}_N^*\setminus\{\emptyset\}}W^*\mathbf{U}^{*w}W\otimes\mathbf{G}^{*w}
\end{eqnarray*}
(note that the sums are finite!). A formal power series $$f(z):=\frac{I_\mathcal{H}}{2}+\sum_{w\in\mathcal{F}_N\setminus\{\emptyset\}}\mathbf{G}^wz^w$$
by Corollary~\ref{cor:coef} belongs to the class $\mathcal{HA}^{\rm nc}_N(\mathcal{H})$. Therefore, for $0<r<1$ one has
\begin{multline*}
I_\mathcal{Y\otimes H}+ \sum_{w\in\mathcal{F}_N\setminus\{\emptyset\}}W^*(r\mathbf{U})^wW\otimes\mathbf{G}^w+
\sum_{w\in\mathcal{F}_N^*\setminus\{\emptyset\}}W^*(r\mathbf{U})^{*w}W\otimes\mathbf{G}^{*w}\\
= (W^*\otimes I_\mathcal{H})2\Re f^{\rm l}(r\mathbf{U})(W\otimes I_\mathcal{H})\geq 0.
\end{multline*}
The sum on the left is finite. Hence, by letting $r\uparrow 1$, we get $p(\mathbf{G})\geq 0$. Thus, $p$ is positive semidefinite on $\mathcal{G}_N$.
\end{proof}
Let us introduce the class ${\rm Nilp}_N$ of $N$-tuples $\mathbf{T}=(T_1,\ldots,T_N)$ of \emph{jointly nilpotent} bounded linear operators on a common Hilbert space, i.e., such that for some $r\in\mathbb{N}$ one has
$$\mathbf{T}^w=0\quad {\rm for\ all}\ w\in\mathcal{F}_N:\ |w|\geq r.$$ 
The minimal such $r$ is called the \emph{rank of joint nilpotency}. Let $\Lambda\subset\mathcal{F}_N$ be an admissible set. We shall say that $\mathbf{T}=(T_1,\ldots,T_N)\in {\rm Nilp}_N$ is an $N$-tuple of \emph{$\Lambda$-jointly nilpotent} operators if 
$$ \mathbf{T}^w=0\quad {\rm for\ all}\ w\in\mathcal{F}_N\setminus\Lambda.$$ 
In this case the rank of joint nilpotency of $\mathbf{T}$ is at most $\max_{w\in\Lambda}|w|+1$. We denote the class of $N$-tuples of $\Lambda$-jointly nilpotent operators by ${\rm Nilp}_N(\Lambda)$. For  $\Lambda_m:=\{ w\in\mathcal{F}_N:\ |w|\leq m\}$, $N$-tuples of $\Lambda_m$-jointly nilpotent operators are exactly those whose rank of joint nilpotency is at most $m+1$. 
\begin{ex}\label{ex:shifts}
Let $\Lambda\subset\mathcal{F}_N$ be an admissible set. Let $\mathcal{H}_\Lambda$ be a finite-dimensional Hilbert space whose orthonormal basis is identified with the set $\Lambda$. For $k=1,\ldots,N$ define the non-commutative backward shifts $S_k\in\mathcal{L}(\mathcal{H}_\Lambda )$ by their action on basis vectors:
$$S_kw =\left\{ \begin{array}{ll}
w' & {\rm if\ } w=g_kw'\ {\rm with\ some\ } w'\in\mathcal{F}_N,\\
0  & {\rm otherwise}.
\end{array}
\right. $$
Since $\Lambda$ is admissible, these operators are correctly defined (if $w=g_kw'$ and $w\in\Lambda$ then $w'\in\Lambda$). The $N$-tuple $\mathbf{S}:=(S_1,\ldots,S_N)$ belongs to the class ${\rm Nilp}_N(\Lambda)$. Indeed, if $v\in\mathcal{F}_N\setminus\Lambda$ and $w\in\Lambda$ then $\mathbf{S}^vw\neq 0$ implies $w=vw'$ with some $w'\in\Lambda$. But in this case (the set $\Lambda$ is admissible!) we get $w\in\mathcal{F}_N\setminus\Lambda$ which is impossible. Thus $\mathbf{S}^v=0$ for every $v\in\mathcal{F}_N\setminus\Lambda$. Since for every $w\in\Lambda$ we have $\mathbf{S}^ww=\emptyset\neq 0$, we obtain that $\mathbf{S}$ doesn't belong to the class ${\rm Nilp}_N(\tilde{\Lambda})$ for any admissible proper subset $\tilde{\Lambda}\subset\Lambda$. We can see also that $\mathbf{S}\in\mathcal{C}^N$: for  arbitrary
$x=\sum_{w\in\Lambda}x_ww\in\mathcal{H}_\Lambda$, with $x_w\in\mathbb{C}$ ($w\in\Lambda$), and $k\in\{ 1,\ldots,N\}$ we get
\begin{eqnarray*}
\| S_kx\|^2 &=& \left\|\sum_{w\in\Lambda}x_wS_kw\right\|^2=\left\|\sum_{w'\in\Lambda:\,g_kw'\in\Lambda}x_{g_kw'}w'\right\|^2\\
&=&
\sum_{w'\in\Lambda:\,g_kw'\in\Lambda}|x_{g_kw'}|^2\leq\sum_{w\in\Lambda}|x_w|^2=\| x\|^2,
\end{eqnarray*}
which means that $S_k$ are contractions.
\end{ex}
\begin{prop}\label{prop:vanish}
Let $\Lambda\subset\mathcal{F}_N$ be an admissible set. If a non-commutative polynomial $p(z)=\sum_{w\in\Lambda}p_wz^w\in\mathcal{L(U,Y)}\left\langle z_1,\ldots,z_N\right\rangle$, with some Hilbert spaces $\mathcal{U}$ and $\mathcal{Y}$, satisfies $p(\mathbf{T})=0$ for an arbitrary $N$-tuple of $\Lambda$-jointly nilpotent contractive $n\times n$ matrices, $n\in\mathbb{N}$,  then $p_w=0$ for all $w\in\Lambda$. Moreover, it suffices to take $n$ equal to $\# (\Lambda)$, the number of words in $\Lambda$.
\end{prop}
\begin{proof}
We have $p(\mathbf{S})=0$, where  $\mathbf{S}$ is the $N$-tuple of backward shifts from Example~\ref{ex:shifts}. Since $\mathbf{S}\in\mathcal{L}(\mathcal{H}_\Lambda)^N$, and $\dim \mathcal{H}_\Lambda =\# (\Lambda)$, one can consider 
$\mathbf{S}$ as an $N$-tuple of $\Lambda$-jointly nilpotent contractive $n\times n$ matrices, with $n=\# (\Lambda)$, as well as $\lambda\mathbf{S}:=(\lambda S_1,\ldots,\lambda S_N)$ for any $\lambda\in\mathbb{D}$. Therefore, a one-variable polynomial
$$r_\mathbf{S}(\lambda)=\sum_{k=0}^mr_{\mathbf{S},k}\lambda^k:=\sum_{k=0}^m\left(\sum_{w\in\Lambda:\,|w|=k}p_w\otimes\mathbf{S}^w\right)\lambda^k=\sum_{w\in\Lambda}p_w\otimes(\lambda\mathbf{S})^w=p(\lambda\mathbf{S}),$$
where $m=\max_{w\in\Lambda}|w|$, vanishes on $\mathbb{D}$, and hence vanishes identically. Thus $$r_{\mathbf{S},k}=\sum_{w\in\Lambda:\,|w|=k}p_w\otimes\mathbf{S}^w=0$$ for $k=0,\ldots,m$.
For a fixed $k$ and any $u\in\mathcal{U}$ and $v\in\Lambda:\ |v|=k$ (the word $v$ is identified with a basis vector in $\mathcal{H}_\Lambda$, or equivalently, with a standard basis vector in $\mathbb{C}^n$), we have:
$$
0=\sum_{w\in\Lambda:\,|w|=k}(p_w\otimes\mathbf{S}^w)(u\otimes v)=\sum_{w\in\Lambda:\,|w|=k}p_wu\otimes\mathbf{S}^wv=p_vu\otimes\emptyset.
$$
Since $\emptyset\neq 0$, we get $p_vu=0$. Since $k\in\{ 0,\ldots,m\}$, $v\in\Lambda$ and $u\in\mathcal{U}$ were chosen arbitrarily, the lemma follows.
\end{proof}
\begin{prop}\label{prop:jn}
Let $m\in\mathbb{N}$. An $N$-tuple $\mathbf{T}\in\mathcal{L(E)}^N$ belongs to the class ${\rm Nilp}_N(\Lambda_m)$ if and only if there exists a decomposition $\mathcal{E}=\mathcal{E}_1\oplus\cdots\oplus\mathcal{E}_{m+1}$ such that the operators $T_k,\ k=1,\ldots,N$, have strictly lower block-triangular form with respect to this decomposition:
$$T_k=\left[\begin{array}{cccc}
0      &   \ldots     &   \ldots     & 0\\
*      & \ddots & \ddots       & \vdots\\
\vdots & \ddots & \ddots & \vdots\\
*      & \ldots & *      & 0 
\end{array}
\right].$$
\end{prop} 
\begin{proof}
Clearly, any $N$-tuple $\mathbf{T}\in\mathcal{L(E)}^N$ of bounded linear operators on a Hilbert space $\mathcal{E}$ which have strictly lower block-triangular form with respect to some decomposition $\mathcal{E}=\mathcal{E}_1\oplus\cdots\oplus\mathcal{E}_{m+1}$ is $\Lambda_m$-jointly nilpotent.

Conversely, let $\mathbf{T}\in\mathcal{L(E)}^N$ be $\Lambda_m$-jointly nilpotent. Set
\begin{eqnarray*}
\mathcal{E}_1 &:=& \left(\bigvee_{w\in\mathcal{F}_N:\,|w|\geq 0}\mathbf{T}^w\mathcal{E}\right)\ominus\left(\bigvee_{w\in\mathcal{F}_N:\,|w|\geq 1}\mathbf{T}^w\mathcal{E}\right)\\
\cdots & \cdots & \cdots \\
\mathcal{E}_{m-1} &:=& \left(\bigvee_{w\in\mathcal{F}_N:\,|w|\geq m-2}\mathbf{T}^w\mathcal{E}\right)\ominus\left(\bigvee_{w\in\mathcal{F}_N:\,|w|\geq m-1}\mathbf{T}^w\mathcal{E}\right)\\
\mathcal{E}_m &:=& \left(\bigvee_{w\in\mathcal{F}_N:\,|w|\geq m-1}\mathbf{T}^w\mathcal{E}\right)\ominus\left(\bigvee_{w\in\mathcal{F}_N:\,|w|=m}\mathbf{T}^w\mathcal{E}\right)\\
\mathcal{E}_{m+1} &:=& \bigvee_{w\in\mathcal{F}_N:\,|w|= m}\mathbf{T}^w\mathcal{E},
\end{eqnarray*}
where $\bigvee_\nu \mathcal{X}_\nu$ denotes the closed linear span of the sets $\mathcal{X}_\nu$ ($\subset\mathcal{E}$).
Then $$\mathcal{E}=\bigoplus_{\nu=1}^{m+1}\mathcal{E}_\nu=\bigvee_{w\in\mathcal{F}_N:\,|w|\geq 0}\mathbf{T}^w\mathcal{E},$$ and $$T_k\mathcal{E}_{j}\subset\bigoplus_{\nu=j+1}^{m+1}\mathcal{E}_\nu=\bigvee_{w\in\mathcal{F}_N:\,|w|\geq j}\mathbf{T}^w\mathcal{E},\quad T_k\mathcal{E}_{m+1}=\{ 0\},\quad k=1,\ldots,N,\ j=1,\ldots,m,$$
which means that $T_k,\ k=1,\ldots,N,$ have strictly lower block-triangular form with respect to the decomposition $\mathcal{E}=\mathcal{E}_1\oplus\cdots\oplus\mathcal{E}_{m+1}$.
\end{proof}
\begin{rem}\label{rem:nm}
The special case of Proposition~\ref{prop:jn} where $\dim\mathcal{E}<\infty$ can be formulated as follows: an $N$-tuple of matrices $\mathbf{T}=(T_1,\ldots,T_N)\in (\mathbf{C}^{n\times n})^N$ is jointly nilpotent, with rank of joint nilpotency at most $m+1$, if and only if $\mathbf{T}$ is unitary similar to an $N$-tuple $\tilde{\mathbf{T}}=(\tilde{T}_1,\ldots,\tilde{T}_N)$ of strictly lower block-triangular $(m+1)\times (m+1)$ matrices,  with not necessarily square blocks. (Here unitary similarity means that there exists a unitary $n\times n$ matrix $U$ such that $T_k=U^{-1}\tilde{T}_kU,\ k=1,\ldots,N$.) This statement is a bit stronger than one in \cite{KV} where only similarity of an $N$-tuple of jointly nilpotent matrices to some $N$-tuple of strictly triangular matrices was mentioned.  
\end{rem}
\begin{lem}\label{lem:dil}
Let $\Lambda\subset\mathcal{F}_N$ be an admissible set, $\mathbf{U}\in\mathcal{U}^N\cap\mathcal{L(K)}^N$, and let $W\in\mathcal{L(Y,K)}$ be an isometry, with Hilbert spaces $\mathcal{Y}$ and $\mathcal{K}$, such that
\begin{equation}\label{cond}
W^*\mathbf{U}^wW=0 \quad {\rm for\ } w\in\mathcal{F}_N\setminus\Lambda.
\end{equation}
Then there exist a Hilbert space $\mathcal{E}$ and an $N$-tuple $\mathbf{T}\in\mathcal{C}^N\cap\mathcal{L(E)}^N$ of $\Lambda$-jointly nilpotent operators such that 
$$
\mathcal{K}\supset\mathcal{E}\supset W\mathcal{Y}
$$
and $\mathbf{U}$ is a unitary dilation of $\mathbf{T}$:
$$\mathbf{T}^w=P_\mathcal{E}\mathbf{U}^w\big|_\mathcal{E},\qquad w\in\mathcal{F}_N.$$
In particular,
\begin{equation}\label{dil}
W^*\mathbf{U}^wW=\tilde{W}^*\mathbf{T}^w\tilde{W},\quad w\in\mathcal{F}_N,
\end{equation}
where $\tilde{W}=P_\mathcal{E}W\in\mathcal{L(Y,E)}$ is an isometry.
If the space $\mathcal{Y}$ is finite-dimensional then one can choose $\mathcal{E}$ finite-dimensional, too.
\end{lem}
\begin{proof}
Define the following subspaces in $\mathcal{K}$:
\begin{eqnarray*}
\mathcal{E}_0 &:=& \bigvee_{w\in\mathcal{F}_N\setminus\Lambda}\mathbf{U}^wW\mathcal{Y},\\
\mathcal{E}   &:=& \left(\bigvee_{w\in\mathcal{F}_N}\mathbf{U}^wW\mathcal{Y}\right)\ominus\left(\bigvee_{w\in\mathcal{F}_N\setminus\Lambda}\mathbf{U}^wW\mathcal{Y}\right),
\end{eqnarray*}
and define the operators $$T_k:=P_\mathcal{E}U_k\big|_\mathcal{E},\quad k=1,\ldots,N.$$
Clearly, $\mathbf{T}=(T_1,\ldots,T_N)\in\mathcal{C}^N$. 
Since both $\mathcal{E}_0$ and $\mathcal{E}_0\oplus\mathcal{E}=\bigvee_{w\in\mathcal{F}_N}\mathbf{U}^wW\mathcal{Y}$ are invariant subspaces for every $U_k,\ k=1,\ldots,N$,  the space $\mathcal{E}$ is an orthogonal difference of two invariant subspaces for these unitary operators, i.e., a \emph{semi-invariant subspace}. Thus, by the Sarason lemma \cite[Lemma~0]{Sar2} $\mathbf{U}$ is a unitary dilation of $\mathbf{T}$. 
Since for every $w\in\mathcal{F}_N\setminus\Lambda$ one has $\mathbf{U}^w\mathcal{E}\subset\mathcal{E}_0$, we get $\mathbf{T}^w=P_\mathcal{E}\mathbf{U}^w\big|_\mathcal{E}=0$, i.e., $\mathbf{T}\in {\rm Nilp}_N(\Lambda)$.
Since by \eqref{cond} $W^*\mathcal{E}_0=\{ 0\}$,  and $W\mathcal{Y}\subset\mathcal{E}_0\oplus\mathcal{E}=\bigvee_{w\in\mathcal{F}_N}\mathbf{U}^wW\mathcal{Y}$, we 
get $W\mathcal{Y}\subset\mathcal{E}$, as desired. Thus, \eqref{dil} holds true as well, with an isometry $\tilde{W}=P_\mathcal{E}W\in\mathcal{L(Y,E)}$. 

In the case where $\dim\mathcal{Y}<\infty$, we may write
$$ \mathcal{E}_0\oplus\mathcal{E}=\left(\bigvee_{w\in\mathcal{F}_N\setminus\Lambda}\mathbf{U}^wW\mathcal{Y}\right)\oplus P_\mathcal{E}\left(\bigvee_{w\in\Lambda}\mathbf{U}^wW\mathcal{Y}\right)=\mathcal{E}_0\oplus P_\mathcal{E}\left(\bigvee_{w\in\Lambda}\mathbf{U}^wW\mathcal{Y}\right),$$ 
and since the set $\Lambda$ is finite, both $\bigvee_{w\in\Lambda}\mathbf{U}^wW\mathcal{Y}$ and $\mathcal{E}=P_\mathcal{E}\left(\bigvee_{w\in\Lambda}\mathbf{U}^wW\mathcal{Y}\right)$ are finite-dimensional subspaces.
\end{proof}
\begin{thm}\label{thm:main}
Problem~\ref{prob:c-gen} has a solution if and only if the polynomial 
\begin{equation}\label{poly}
p(z):=\frac{c_\emptyset}{2}+\sum_{w\in\Lambda\setminus\{\emptyset\} }c_wz^w
\end{equation}
satisfies $\Re p(\mathbf{T})\geq 0$ for every $N$-tuple $\mathbf{T}$ of $\Lambda$-jointly nilpotent contractive operators. Moreover, for the solvability of Problem~\ref{prob:c-gen} it is enough to assume that $\Re p(\mathbf{T})\geq 0$ holds for every $N$-tuple $\mathbf{T}$ of $\Lambda$-jointly nilpotent contractive $n\times n$ matrices, for all $n\in\mathbb{N}$.
\end{thm}
\begin{proof}
If Problem~\ref{prob:c-gen} has a solution $f\in\mathcal{HA}_N^{\rm nc}(\mathcal{Y})$ then  for any $\mathbf{C}\in\mathcal{D}^N\cap\mathcal{L(E)}^N$, with a Hilbert space $\mathcal{E}$, the series $$f(\mathbf{C}):=\sum_{w\in\mathcal{F}_N}f_w\otimes\mathbf{C}^w$$
converges in the operator norm, and $\Re f(\mathbf{C})\geq 0$. If $\mathbf{T}\in\mathcal{C}^N\cap\mathcal{L(E)}^N$ is an $N$-tuple of $\Lambda$-jointly nilpotent operators then so is $r\mathbf{T}=(rT_1,\ldots,rT_N)\in\mathcal{D}^N\cap\mathcal{L(E)}^N$ for every $r:\ 0<r<1$. Therefore, $$\Re p(r\mathbf{T})=\Re f(r\mathbf{T})=\Re\left(\frac{c_\emptyset\otimes I_\mathcal{E}}{2}+\sum_{w\in\Lambda\setminus\{ \emptyset\}}c_w\otimes(r\mathbf{T})^w\right)\geq 0.$$
By letting $r\uparrow 1$, we obtain $\Re p(\mathbf{T})\geq 0$.

For the converse direction, let us consider first the case where $c_\emptyset=I_\mathcal{Y}$. Let $\Re p(\mathbf{T})\geq 0$ hold for every $N$-tuple $\mathbf{T}$ of $\Lambda$-jointly nilpotent contractive $n\times n$ matrices, for all $n\in\mathbb{N}$. Let $\mathcal{A}_{\mathcal{G}_N}$ be the $C^*$-algebra obtained as  the norm completion of the quotient of unital $*$-algebra $\mathcal{A}_N$ (which has been introduced in the proof of Lemma~\ref{lem:main} above) with the seminorm
$$\| q\| :=\sup_{\mathbf{G}\in\mathcal{G}_N}\| q(\mathbf{G})\|=\sup_{\mathbf{G}\in\mathcal{G}^N}\| q(G_1,\ldots,G_N,G_1^*,\ldots,G_N^*)\|,$$
by the two-sided ideal of elements of zero seminorm. 

Let us show that the restriction of the quotient map above to the subspace $\mathcal{B}_N\subset\mathcal{A}_N$ of polynomials of the form \eqref{q}
 is injective, i.e., that if $q\in\mathcal{B}_N$ is non-zero then the corresponding coset $[q]\in\mathcal{A}_{\mathcal{G}_N}$ is non-zero. Indeed, if $[q]=[0]$ then $q(\mathbf{G})=0$ for every $\mathbf{G}\in\mathcal{G}_N\cap\mathcal{L}(\mathcal{H})^N$, with a Hilbert space $\mathcal{H}$. Define $\tilde{q}(\hat{z}):=1+q(\hat{z})$. Then $\tilde{q}(\mathbf{G})=I_{\mathcal{H}}>0$ for every $\mathbf{G}\in\mathcal{G}_N\cap\mathcal{L}(\mathcal{H})^N$. In particular, $\tilde{q}(\mathbf{G}_{\zeta})=1$ for $\mathbf{G}_{\zeta}:=(\zeta,0,\ldots,0)\in\mathcal{G}_N\cap\mathbb{C}^N,\ \zeta\in\mathbb{T}$. Therefore,
$$\tilde{q}_\emptyset=\int_{\mathbb{T}}\tilde{q}(\mathbf{G}_\zeta)\,d\zeta =1,$$
and $q_\emptyset=0$.
 By Part~II of Lemma~\ref{lem:main}, there exist a Hilbert space $\mathcal{K}$, an $N$-tuple of operators $\mathbf{U}\in\mathcal{U}^N\cap\mathcal{L}(\mathcal{K})^N$, and an isometry $W\in\mathcal{L}(\mathbb{C},\mathcal{K})$ such that 
	\begin{eqnarray*}
	\tilde{q}_w=\tilde{q}_{w^*}^* &=&  W^*\mathbf{U}^wW,\qquad w\in\mathcal{F}_N:\,|w|\leq m,\\
	0 &=& W^*\mathbf{U}^wW,\qquad w\in \mathcal{F}_N:\, |w|>m.
\end{eqnarray*}
Then $$q_w=\tilde{q}_w=\tilde{q}_{w^*}^*=q_{w^*}^*=W^*\mathbf{U}^w\tilde{W},\qquad w\in\mathcal{F}_N:\,0<|w|\leq m.$$
For an arbitrary $\mathbf{G}\in\mathcal{G}_N\cap\mathcal{L}(\mathcal{H})^N$ the one-variable trigonometric polynomial $$t_\mathbf{G}(\zeta):=q(\zeta\mathbf{G})=\sum_{k=1}^m\left(\sum_{w\in\mathcal{F}_N:\,|w|=k}q_w\mathbf{G}^w\right)\zeta^k+\sum_{k=1}^m\left(\sum_{w\in\mathcal{F}_N^*:\,|w|=k}q_w\mathbf{G}^{*w}\right)\bar{\zeta}^k$$
is identically zero, which implies in particular
$$\sum_{w\in\mathcal{F}_N:\,|w|=k}q_w\mathbf{G}^w=0,\qquad  k=1,\ldots,m.$$
For arbitrary $n\in\mathbb{N}$ and  $\tilde{\mathbf{U}}\in\mathcal{U}^N\cap {(\mathbb{C}^{n\times n})}^N$, by virtue of Proposition~\ref{prop:schur} we have $\mathbf{G}\stackrel{\circ}{\otimes }\tilde{\mathbf{U}}\in\mathcal{G}_N$, therefore
$$\sum_{w\in\mathcal{F}_N:\,|w|=k}q_w(\mathbf{G}\stackrel{\circ}{\otimes}\tilde{\mathbf{U}})^w=\sum_{w\in\mathcal{F}_N:\,|w|=k}q_w\mathbf{G}^w\otimes\tilde{\mathbf{U}}^w=0,\qquad  k=1,\ldots,m.$$
Since $\mathcal{U}^N\cap {(\mathbb{C}^{n\times n})}^N$ is a uniqueness set for holomorphic functions of matrix entries (see e.g., \cite{Sh}), the non-commutative polynomial $\sum_{w\in\mathcal{F}_N:\,|w|=k}q_w\mathbf{G}^wz^w$ vanishes on the all of ${(\mathbb{C}^{n\times n})}^N$, for every $k\in\{ 1,\ldots,m\}$ and $n\in\mathbb{N}$:
$$\sum_{w\in\mathcal{F}_N:\,|w|=k}q_w\mathbf{G}^w\otimes \mathbf{Z}^w=0,\qquad \mathbf{Z}\in {(\mathbb{C}^{n\times n})}^N.$$ Thus, by the Amitsur--Levitzki theorem (see \cite[pp. 22--23]{Ro}), $$q_w\mathbf{G}^w=0,\qquad w\in\mathcal{F}_N:\ 0<|w|\leq m.$$ For every $w\in\mathcal{F}_N:\ 0<|w|\leq m$ one can find $\mathbf{G}\in\mathcal{G}_N$ such that 
$\mathbf{G}^w\neq 0$. Indeed, the non-commutative polynomial $\psi(z)=1+z^w$ belongs to the class $\mathcal{HA}_N^{{\rm nc}, 1}=\mathcal{HA}_N^{{\rm nc},1}(\mathbb{C})$, thus by Corollary~\ref{cor:coef}, there exist a Hilbert space $\mathcal{H}$, an $N$-tuple of operators  $\mathbf{G}\in\mathcal{G}_N\cap\mathcal{L}(\mathcal{H})^N$, and an isometry $V\in\mathcal{L}(\mathbb{C}, \mathcal{H})$ such that $\psi_w=1=2V^*\mathbf{G}^wV$, which implies $\mathbf{G}^w\neq 0$. Therefore, $q_w=0$ for all $w\in\mathcal{F}_N:\ 0<|w|\leq m$. Since we have shown already that $q_\emptyset=0$, and $q_w=q_{w^*}^*$ for $w\in\mathcal{F}_N:\ 0<|w|\leq m$, we get $q=0$.

Denote by $\mathcal{B}_\Lambda\subset\mathcal{B}_N$ the finite-dimensional subspace of polynomials of the form
\begin{equation}\label{ql}
q(\hat{z})=q_\emptyset+\sum_{w\in\Lambda\setminus\{\emptyset\} }q_wz^w+\sum_{w\in\Lambda^*\setminus\{\emptyset\} }q_wz^{*w},
\end{equation}
and let $\mathcal{B}_{\Lambda,\mathcal{G}_N}$ be the image of the subspace $\mathcal{B}_\Lambda$ under the quotient map above. This subspace of the $C^*$-algebra $\mathcal{A}_{\mathcal{G}_N}$ is selfadjoint, i.e., $\mathcal{B}_{\Lambda,\mathcal{G}_N}^*=\mathcal{B}_{\Lambda,\mathcal{G}_N}.$ Define the linear map $\varphi:\ \mathcal{B}_{\Lambda,\mathcal{G}_N}\to\mathcal{L(Y)}$ by $\varphi([1])=I_\mathcal{Y}$, $\varphi([z^w])=c_w$ for $w\in\Lambda\setminus\{ \emptyset\}$, and $\varphi([z^{*w}])=c_{w^*}^*$ for $w\in\Lambda^*\setminus\{ \emptyset\}$. By the previous paragraph together with the above mentioned Amitsur--Levitzki theorem, this linear map is correctly defined. Let us show that $\varphi$ is completely positive. Let $n\in\mathbb{N}$, and let $[q]\in\mathbb{C}^{n\times n}\otimes\mathcal{B}_{\Lambda,\mathcal{G}_N}$ be a positive element of the $C^*$-algebra $\mathbb{C}^{n\times n}\otimes\mathcal{A}_{\mathcal{G}_N}$, i.e.,
$[q]=[h]^*[h]$ with some $[h]\in\mathbb{C}^{n\times n}\otimes\mathcal{A}_{\mathcal{G}_N}$. One can think of $[q]$ as of the $n\times n$ matrix $([q]_{ij})_{i,j=1,\ldots,n}$ whose entries $[q]_{ij}=[q_{ij}]\in\mathcal{B}_{\Lambda,\mathcal{G}_N}$ and $q_{ij}\in\mathcal{B}_\Lambda$, and thus $q\in\mathbb{C}^{n\times n}\otimes\mathcal{B}_\Lambda$ is a polynomial of the form \eqref{ql} with the coefficients from $\mathbb{C}^{n\times n}$. Let us observe that by virtue of the definition of the $C^*$-algebra $\mathcal{A}_{\mathcal{G}_N}$, for an arbitrary $[x]\in\mathcal{A}_{\mathcal{G}_N}$ its values on $\mathcal{G}_N$ are well defined. In particular, if $x\in\mathcal{B}_\Lambda$ then $[x](\mathbf{G})=x(\mathbf{G})$ for any $\mathbf{G}\in\mathcal{G}_N$. Therefore, for an arbitrary $[x]=([x]_{ij})_{i,j=1,\ldots,n}=([x_{ij}])_{i,j=1,\ldots,n}\in\mathbb{C}^{n\times n}\otimes\mathcal{A}_{\mathcal{G}_N}$ one defines correctly $[x](\mathbf{G}):=([x_{ij}](\mathbf{G}))_{i,j=1,\ldots,n},\ \mathbf{G}\in\mathcal{G}_N$. In particular, if $x=(x_{ij})_{i,j=1,\ldots,n}\in\mathbb{C}^{n\times n}\otimes\mathcal{B}_\Lambda$ then $[x](\mathbf{G})=x(\mathbf{G})=(x_{ij}(\mathbf{G}))_{i,j=1,\ldots,n}$ for any $\mathbf{G}\in\mathcal{G}_N$. Since $q(\mathbf{G})=[q](\mathbf{G})=[h](\mathbf{G})^*[h](\mathbf{G})$ is positive semidefinite for every  $\mathbf{G}\in\mathcal{G}_N$, it follows that the polynomial $q$ (with the coefficients from $\mathbb{C}^{n\times n}$) is positive semidefinite on $\mathcal{G}_N$. In particular, $q(\mathbf{G}_{\zeta})\geq 0$ for $\mathbf{G}_{\zeta}:=(\zeta,0,\ldots,0)\in\mathcal{G}_N\cap\mathbb{C}^N,\ \zeta\in\mathbb{T}$. Therefore,
$$q_\emptyset=\int_{\mathbb{T}}q(\mathbf{G}_\zeta)\,d\zeta\geq 0.$$
  If $q_\emptyset=I_n$, then by Part II of Lemma~\ref{lem:main} there exist a Hilbert space $\mathcal{K}$, an $N$-tuple of operators $\mathbf{U}\in\mathcal{U}^N\cap\mathcal{L}(\mathcal{K})^N$, and an isometry $W\in\mathcal{L}(\mathbb{C}^n,\mathcal{K})$ such that 
\begin{eqnarray*}
q_w=q_{w^*}^* &=& W^*\mathbf{U}^wW,\qquad w\in\Lambda,\\
0   &=& W^*\mathbf{U}^wW,\qquad w\in\mathcal{F}_N\setminus\Lambda.
\end{eqnarray*}
Then we have
\allowdisplaybreaks
\begin{eqnarray*}
\varphi_n([q]) &=& ({\rm id}_n\otimes \varphi)\left(I_n\otimes [1]+\sum_{w\in\Lambda\setminus\{\emptyset\} }q_w\otimes [z^w]+ \sum_{w\in\Lambda^*\setminus\{\emptyset\} }q_w\otimes [z^{*w}] \right)\\ 
&=& I_{\mathbb{C}^n\otimes \mathcal{Y}}+\sum_{w\in\Lambda\setminus\{\emptyset\} }W^*\mathbf{U}^wW\otimes c_w
+ \sum_{w\in\Lambda^*\setminus\{\emptyset\} }W^*\mathbf{U}^{*w}W\otimes c_{w^*}^*\\
&=& (W^*\otimes I_\mathcal{Y})2\Re p^{\rm l}(\mathbf{U})(W\otimes I_\mathcal{Y})
\end{eqnarray*}
By Lemma~\ref{lem:dil}, there exist a finite-dimensional Hilbert space $\mathcal{E}$ and an $N$-tuple $\mathbf{T}\in\mathcal{C}^N\cap\mathcal{L(E)}^N$ of $\Lambda$-jointly nilpotent operators such that 
\eqref{dil} holds
with an isometry $\tilde{W}\in\mathcal{L}(\mathbb{C}^n,\mathcal{E})$. Thus,
$$\varphi_n([q])=(\tilde{W}^*\otimes I_\mathcal{Y})2\Re p^{\rm l}(\mathbf{T})(\tilde{W}\otimes I_\mathcal{Y})$$
is a positive semidefinite operator by the assumption that $\Re p(\mathbf{T})\geq 0$ or, equivalently, $\Re 
p^{\rm l}(\mathbf{T})\geq 0$.
In the case where $q_\emptyset >0$ we can define $\tilde{q}(\hat{z}):=q_\emptyset^{-1/2}q(\hat{z})q_\emptyset^{-1/2}$. Since $\varphi_n([\tilde{q}])\geq 0$, we get $$\varphi_n([q])=(q_\emptyset^{1/2}\otimes I_\mathcal{Y})\varphi_n([\tilde{q}])(q_\emptyset^{1/2}\otimes I_\mathcal{Y})\geq 0.$$
In the case where the matrix $q_\emptyset$ is degenerate we set $q_\epsilon(\hat{z}) :=\epsilon I_n+q(\hat{z})$ for  $\epsilon >0$. Then $q_\epsilon$ is positive definite on $\mathcal{G}_N$ and ${(q_\epsilon)}_\emptyset=\epsilon I_n+q_\emptyset >0$. Since $\varphi_n([q_\epsilon])\geq 0$, we get
 $$
\varphi_n([q]) = \lim_{\epsilon\downarrow 0}\varphi_n([q_\epsilon ])\geq 0.
$$
Finally, we have obtained that $\varphi: \mathcal{B}_{\Lambda,\mathcal{G}_N}\rightarrow \mathcal{L(Y)}$ is completely positive.

 Since we have $\varphi ([1])=I_\mathcal{Y}$, by the Arveson extension theorem \cite{Arv1} there exists a completely positive map $\widetilde{\varphi}:\mathcal{A}_{\mathcal{G}_N}\rightarrow \mathcal{L(Y)}$ which extends $\varphi$. By the Stinespring theorem \cite{St}, there exists a $*$-representation $\pi$ of $\mathcal{A}_{\mathcal{G}_N}$ in some Hilbert space $\mathcal{H}$ and an isometry $V\in\mathcal{L(Y,H)}$ such that
$$\widetilde{\varphi}(a)=V^*\pi(a)V,\qquad a\in\mathcal{A}_{\mathcal{G}_N}.$$
In particular, we get
\begin{equation}
c_w = \varphi([z^w])=\widetilde{\varphi}([z^w])=V^*\pi([z^w])V=V^*\mathbf{G}^wV,\quad w\in\Lambda,\label{II1''}\\
\end{equation}
where we set $\mathbf{G}:=(\pi([z_1]),\ldots,\pi([z_N]))$. We have $\mathbf{G}\in\mathcal{G}_N$. Indeed, since by Proposition~\ref{prop:gn} we have $$\left\| \left[1-\sum_{k=1}^Nz_k^*z_k\right]\right\|_{\mathcal{A}_{\mathcal{G}_N}}=\left\| \left[1-\sum_{k=1}^Nz_kz_k^*\right]\right\|_{\mathcal{A}_{\mathcal{G}_N}}=0,$$
and $$\|[z_k^*z_j]\|_{\mathcal{A}_{\mathcal{G}_N}}=\|[z_kz_j^*]\|_{\mathcal{A}_{\mathcal{G}_N}}=0,\quad k\neq j,$$
 we get $[\sum_{k=1}^Nz_k^*z_k]=[\sum_{k=1}^Nz_kz_k^*]=[1]$, and $[z_k^*z_j]=[z_kz_j^*]=[0]$ for $k\neq j$. Hence
\begin{eqnarray*}
	\sum_{k=1}^NG_k^*G_k &=& \sum_{k=1}^N\pi([z_k])^*\pi([z_k])=\pi\left(\left[ \sum_{k=1}^Nz_k^*z_k\right]\right)=\pi([1])=I_\mathcal{H},\\
	G_k^*G_j &=& \pi([z_k])^*\pi([z_j])=\pi([z_k^*z_j])=\pi([0])=0,\qquad k\neq j,\\
	\sum_{k=1}^NG_kG_k^* &=& \sum_{k=1}^N\pi([z_k])\pi([z_k])^*=\pi\left(\left[ \sum_{k=1}^Nz_kz_k^*\right]\right)=\pi([1])=I_\mathcal{H},\\
	G_kG_j^* &=& \pi([z_k])\pi([z_j])^*=\pi([z_kz_j^*])=\pi([0])=0,\qquad k\neq j,
\end{eqnarray*} 
which means that $\mathbf{G}\in\mathcal{G}_N$, according to Proposition~\ref{prop:gn}.
	Finally, since \eqref{II1''} coincides with \eqref{mom-nc}, from Theorem~\ref{thm:g} we obtain that  Problem~\ref{prob:c} has a solution.

Consider now the case where $c_\emptyset >0$. Set $\widetilde{c_w}:=c_\emptyset^{-1/2}c_wc_\emptyset^{-1/2},\ w\in\Lambda$, and $\widetilde{p}(z):=c_\emptyset^{-1/2}p(z)c_\emptyset^{-1/2}$, where $p(z)$ is given by \eqref{poly}. Clearly, $\widetilde{c_\emptyset }=I_\mathcal{Y}$. Since Problem~\ref{prob:c} for the data $\widetilde{c_w},\ w\in\Lambda$, is solvable if and only if $\Re\widetilde{p}(\mathbf{T})\geq 0$ for every $N$-tuple $\mathbf{T}$ of contractive $\Lambda$-jointly nilpotent square matrices of same size, and since $\Re\widetilde{p}(\mathbf{T})\geq 0\Leftrightarrow\Re p(\mathbf{T})\geq 0$ and $f\in\mathcal{HA}_N^{\rm nc}(\mathcal{Y})\Leftrightarrow c_\emptyset^{-1/2}fc_\emptyset^{-1/2}\in\mathcal{HA}_N^{\rm nc}(\mathcal{Y})$, Problem~\ref{prob:c-gen} for the data $c_\emptyset>0,\ c_w\ (w\in\Lambda\setminus\{\emptyset\} )$ is solvable if and only if $\Re p(\mathbf{T})\geq 0$ for every $N$-tuple $\mathbf{T}$ of contractive $\Lambda$-jointly nilpotent square matrices of same size.

Consider now the general case $c_\emptyset\geq 0$. Suppose that $\Re p(\mathbf{T})\geq 0$ for every $N$-tuple $\mathbf{T}$ of contractive $\Lambda$-jointly nilpotent $n\times n$ matrices, for all $n\in\mathbb{N}$. Then for any such $\mathbf{T}$ the polynomial $$f_\mathbf{T}(\lambda):=p(\lambda\mathbf{T})=\frac{c_\emptyset\otimes I_n}{2}+\sum_{k=1}^m \left(\sum_{w\in\Lambda:\,|w|=k}c_w\otimes\mathbf{T}^w\right)\lambda^k,\quad \lambda\in\mathbb{D},$$
where $m=\max_{w\in\Lambda}|w|$,
belongs to the Herglotz class $\mathcal{H}_1(\mathcal{Y}\otimes\mathbb{C}^{n})$. Since its coefficients are   $\frac{c_0(\mathbf{T})}{2}:=(f_\mathbf{T})_0=\frac{c_\emptyset\otimes I_n}{2}$, $c_k(\mathbf{T}):=(f_\mathbf{T})_k=\sum_{w\in\Lambda:\,|w|=k}c_w\otimes\mathbf{T}^w,\ k=1,\ldots,m,$  from the Carath\'{e}odory--Toeplitz criterion of solvability of the one-variable Carath\'{e}odory problem with data $c=\{ c_k(\mathbf{T})\}_{k=0,\ldots,m}$ we obtain $T_c\geq 0$, where the operator block matrix $T_c$ is defined by \eqref{t}.
The condition $T_c\geq 0$ implies for $k=1,\ldots,m$ the following inequalities:
\begin{equation}\label{ineq}
|\left\langle c_k(\mathbf{T})x,y\right\rangle |^2\leq\left\langle (c_\emptyset\otimes I_n) x,x\right\rangle\left\langle (c_\emptyset\otimes I_n) y,y\right\rangle,\quad x,y\in\mathcal{Y}\otimes\mathbb{C}^n,
\end{equation}
which yield  $\ker c_\emptyset\otimes\mathbb{C}^n\subset \ker c_k(\mathbf{T})$, and $\ker c_\emptyset\otimes\mathbb{C}^n\subset\ker c_k(\mathbf{T})^*$. Let $x\in\ker c_\emptyset$ and $k\in\{ 1,\ldots,m\}$. Then the non-commutative polynomial $\sum_{w\in\Lambda:\,|w|=k}(c_wx)z^w$ with coefficients in $\mathcal{Y}\cong\mathcal{L}(\mathbb{C},\mathcal{Y})$ vanishes on  $N$-tuples of contractive  $\Lambda$-jointly nilpotent $n\times n$ matrices, for every $n\in\mathbb{N}$. By Proposition~\ref{prop:vanish}, $c_wx=0$ for all $w\in\Lambda:\ |w|=k$.
 Therefore, for every $w\in\Lambda$ we have $\ker c_\emptyset\subset\ker c_w$. Analogously,  $\ker c_\emptyset\subset\ker c_w^*$ for every $w\in\Lambda$. We obtain that our data of Problem~\ref{prob:c-gen} have the following operator block matrix form with respect to the decomposition $\mathcal{Y}=\ker c_\emptyset\otimes {\rm ran}\,c_\emptyset$:
 $$c_w=\left[\begin{array}{cc}
0 & 0\\
0 & c_w^{(22)}\end{array}\right],\quad w\in\Lambda.$$
 Since $c_\emptyset^{(22)}>0$, and the polynomial $$p^{(22)}(z):=\frac{c_\emptyset^{(22)}}{2}+\sum_{w\in\Lambda\setminus\{\emptyset\} }c_w^{(22)}z^w$$
 satisfies the condition that $\Re p^{(22)}(\mathbf{T})\geq 0$ for every $N$-tuple $\mathbf{T}$ of contractive $\Lambda$-jointly nilpotent square matrices of same size, by the result of the previous paragraph, Problem~\ref{prob:c-gen} for the data $c_w^{(22)},\ w\in\Lambda$, has a solution $f^{(22)}\in\mathcal{HA}_N^{\rm nc}({\rm ran}\,c_\emptyset)$. Then Problem~\ref{prob:c-gen} for the data $c_w,\ w\in\Lambda$, has a solution
 $$f=\left[\begin{array}{cc}
0 & 0\\
0 & f^{(22)}\end{array}\right]\in\mathcal{HA}_N^{\rm nc}(\mathcal{Y}).$$
\end{proof}

Let us remark that the condition that $\Re p(\mathbf{T})\geq 0$ for every $N$-tuple of $\tilde{\Lambda}$-jointly nilpotent contractive $n\times n$ matrices, for all $n\in\mathbb{N}$, where $\tilde{\Lambda}\supset\Lambda$ is a wider admissible set, is sufficient for the solvability of Problem~\ref{prob:c-gen}. For instance, one might find convenient to test this condition for the set $\tilde{\Lambda}=\Lambda_m$, with $m=\max_{w\in\Lambda}|w|$, and use the structure of $\Lambda_m$-jointly nilpotent matrices described in Remark~\ref{rem:nm}. However, one should remember that in general this condition is not necessary for the solvability of Problem~\ref{prob:c-gen}.
\begin{ex}\label{ex:1}
Let $\Lambda=\{\emptyset,g_1,g_2,g_1g_2,g_2g_1\}\in\mathcal{F}_2$, and the scalar data of the Carath\'{e}odory problem are
$c_\emptyset=1,\ c_{g_1}=c_{g_2}=\frac{1}{2},\ c_{g_1g_2}=c_{g_2g_1}=\frac{1}{4}$. Then the formal power series
$$\frac{1}{2}\left(1+\frac{z_1+z_2}{2}\right)\left(1-\frac{z_1+z_2}{2}\right)^{-1}=\frac{1}{2}+\sum_{j=1}^\infty\left(\frac{z_1+z_2}{2}\right)^j\in\mathcal{HA}_2^{\rm nc}$$
is a solution to this problem. By Theorem~\ref{thm:main}, for the polynomial $$p(z):=\frac{1}{2}+\frac{z_1+z_2}{2}+\frac{z_1z_2+z_2z_1}{4}$$ and  for every pair of $n\times n$ matrices $\mathbf{T}=(T_1,T_2)\in\mathcal{C}^2\cap {\rm Nilp}_2(\Lambda),\ n\in\mathbb{N}$, one has $\Re p(\mathbf{T})\geq 0$.
Let $\tilde{\Lambda}:=\{\emptyset,g_1,g_2,g_1g_2,g_2g_1,g_1^2\}\supset\Lambda$, and $\mathbf{T}=(T_1,T_2)\in (\mathbb{C}^{3\times 3})^2$ be given by
$$T_1:=\left[\begin{array}{ccc}
0 & 1 & 0\\
0 & 0 & 1\\
0 & 0 & 0\end{array}\right],\quad T_2:=\left[\begin{array}{ccc}
0 & 1 & 0\\
0 & 0 & 0\\
0 & 0 & 0\end{array}\right].$$
It is easy to see that $\mathbf{T}\in\mathcal{C}^2\cap {\rm Nilp}_2(\tilde{\Lambda})$, and $\mathbf{T}\notin {\rm Nilp}_2(\Lambda)$ since $T_1^2\neq 0$. Since $$\det (2\Re p(\mathbf{T}))=\left|\begin{array}{ccc}
1 & 1 & 1/4\\
1 & 1 & 1/2\\
1/4 & 1/2 & 1\end{array}\right| =-\frac{1}{16}<0,$$
the condition $\Re p(\mathbf{T})\geq 0$ is not fulfilled. The same is true for the (admissible) set $\tilde{\Lambda}:=\{\emptyset,g_1,g_2,g_1g_2,g_2g_1,g_2^2\}\supset\Lambda$, where we take the previous example of $\mathbf{T}=(T_1,T_2)$ and interchange $T_1\leftrightarrow T_2$. Clearly, the condition $\Re p(\mathbf{T})\geq 0$ can not be fulfilled for all pairs of $n\times n$ matrices $\mathbf{T}\in\mathcal{C}^2\cap {\rm Nilp}_2(\Lambda_2),\ n\in\mathbb{N}$.
\end{ex}
The following example shows that sometimes the above mentioned condition for a wider set $\tilde{\Lambda}\supset \Lambda$
is necessary for the solvability of Problem~\ref{prob:c-gen}.
\begin{ex}\label{ex:2}
Let $\Lambda=\{\emptyset,g_1,g_1^2,\ldots,g_1^m\}\subset\mathcal{F}_N$ and $c_\emptyset\geq 0,c_{g_1},c_{g_1^2},\ldots,c_{g_1^m}\in\mathcal{L(Y)}$, with some $m\in\mathbb{N}$ and some Hilbert space $\mathcal{Y}$. Then the class ${\rm Nilp}_N(\Lambda)$ consists of $N$-tuples of operators of the form $\mathbf{T}=(T_1,0,\ldots,0)$, where $T_1$ is a nilpotent operator. The condition that $\Re p(\mathbf{T})=\Re (\frac{c_\emptyset\otimes I}{2}+\sum_{j=1}^mc_{g_1^j}\otimes T_1^j)\geq 0$ for every nilpotent contractive square matrix $T_1$ is necessary and sufficient for the solvability of Problem~\ref{prob:c-gen} for these data (and equivalent to the Carath\'{e}odory--Toeplitz criterion for the solvability of a one-variable Carath\'{e}odory problem, see Section~\ref{sec:intro}). Then for every admissible set $\tilde{\Lambda}\supset\Lambda$ such that $g_1^{m+1}\notin\tilde{\Lambda}$ the condition $\Re p(\mathbf{T})\geq 0$ is fulfilled for every $N$-tuple of matrices $\mathbf{T}\in\mathcal{C}^N\cap {\rm Nilp}_N(\tilde{\Lambda})$. In particular it is fulfilled for $\tilde{\Lambda}=\Lambda_m$.
\end{ex}
Thus, natural open questions are the following.
For which admissible sets $\Lambda\subset\mathcal{F}_N$ and data $c_w,\ w\in\Lambda$, the condition that $\Re p(\mathbf{T})\geq 0$ for every $N$-tuple $\mathbf{T}$ of $n\times n$ matrices $\mathbf{T}\in\mathcal{C}^N\cap {\rm Nilp}_N(\Lambda_m),\ n\in\mathbb{N}$, where $m=\max_{w\in\Lambda}|w|$, is necessary for the solvability of Problem~\ref{prob:c-gen}?  Which admissible sets $\Lambda\subset\mathcal{F}_N$ are maximal in the sense that, for a certain choice of problem data $c_w,\ w\in\Lambda$, the condition above fails not only for $\Lambda_m\supset\Lambda$ but also for every admissible set $\tilde{\Lambda}\supset\Lambda$?

\section{The Carath\'{e}odory--Fej\'{e}r interpolation problem}\label{sec:prob-cf}
Recall that the \emph{non-commutative Schur--Agler class} $\mathcal{SA}_N^{\rm nc}(\mathcal{U,Y})$ consists of formal power series $F(z)=\sum_{w\in\mathcal{F}_N}F_wz^w\in\mathcal{L(U,Y)}\left\langle \left\langle z_1,\ldots,z_N\right\rangle\right\rangle$ such that  for every $\mathbf{C}\in\mathcal{D}^N$ (or equivalently, for every $\mathbf{C}\in\mathcal{D}^N_{\rm matr}$, see \cite{AK2})  the series $\sum_{w\in\mathcal{F}_N}F_w\otimes\mathbf{C}^w$ converges in the operator norm to the contractive operator $F(\mathbf{C})$.

Let us pose now the \emph{Carath\'{e}odory--Fej\'{e}r  problem in the class $\mathcal{SA}_N^{\rm nc}(\mathcal{U,Y})$}.
\begin{prob}\label{prob:cf}
Let $\Lambda\subset\mathcal{F}_N$ be an admissible set. Given a collection of operators $\{ s_w\}_{w\in\Lambda}\in\mathcal{L(U,Y)}$, find 
$F\in\mathcal{SA}_N^{\rm nc}(\mathcal{U,Y})$ such that
$$
F_w=s_w,\quad w\in\Lambda.
$$
\end{prob}
\begin{thm}\label{thm:cf}
Problem~\ref{prob:cf} has a solution if and only if the polynomial 
\begin{equation}\label{q-poly}
q(z):=\sum_{w\in\Lambda}s_wz^w
\end{equation}
satisfies $\| q(\mathbf{T})\|\leq 1$ for every $N$-tuple $\mathbf{T}$ of $\Lambda$-jointly nilpotent contractive operators. Moreover, for the solvability of Problem~\ref{prob:cf} it is enough to assume that $\| q(\mathbf{T})\|\leq 1$ holds for every $N$-tuple $\mathbf{T}$ of $\Lambda$-jointly nilpotent contractive $n\times n$ matrices, for all $n\in\mathbb{N}$.
\end{thm}
\begin{proof}
If Problem~\ref{prob:cf} has a solution $F\in\mathcal{SA}_N^{\rm nc}(\mathcal{Y})$ then  for any $\mathbf{C}\in\mathcal{D}^N\cap\mathcal{L(E)}^N$, with a Hilbert space $\mathcal{E}$, the series $F(\mathbf{C}):=\sum_{w\in\mathcal{F}_N}F_w\otimes\mathbf{C}^w$
converges in the operator norm, and $\| F(\mathbf{C})\|\leq 1$. If $\mathbf{T}\in\mathcal{C}^N\cap\mathcal{L(E)}^N$ is an $N$-tuple of $\Lambda$-jointly nilpotent operators then so is $r\mathbf{T}=(rT_1,\ldots,rT_N)\in\mathcal{D}^N\cap\mathcal{L(E)}^N$ for every $r:\ 0<r<1$. Therefore, $$\| q(r\mathbf{T})\|=\| F(r\mathbf{T})\|=\left\|\sum_{w\in\Lambda}s_w\otimes(r\mathbf{T})^w\right\|\leq 1.$$
By letting $r\uparrow 1$, we obtain $\| q(\mathbf{T})\|\leq 1$.

For the converse direction, let us consider first the case where $\mathcal{U}=\mathcal{Y}$ and $-I_\mathcal{Y}\leq s_\emptyset=s_\emptyset^*\leq 0$. Then the operator $I_\mathcal{Y}-s_\emptyset$ is boundedly invertible, and $(I_\mathcal{Y}+s_\emptyset)(I_\mathcal{Y}-s_\emptyset)^{-1}\geq 0$. Moreover, a formal power series $h(z):=(I_\mathcal{Y}+q(z))(I_\mathcal{Y}-q(z))^{-1}\in\mathcal{L(Y)}\left\langle \left\langle z_1,\ldots,z_N\right\rangle\right\rangle$ is well defined. Suppose that $\| q(\mathbf{T})\|\leq 1$ holds for every $N$-tuple $\mathbf{T}$ of $\Lambda$-jointly nilpotent contractive $n\times n$ matrices, for all $n\in\mathbb{N}$. Then
$\Re h(\mathbf{T})\geq 0$ for such a $\mathbf{T}$ (here $ h(\mathbf{T})=\sum_{w\in\Lambda}h_w\otimes\mathbf{T}^w$ is well defined). Define $c_\emptyset:=2h_\emptyset=2(I_\mathcal{Y}+s_\emptyset)(I_\mathcal{Y}-s_\emptyset)^{-1}\geq 0$, $c_w:=h_w$ for $w\in\Lambda\setminus\{\emptyset\}$. For the polynomial $p(z)$ defined by \eqref{poly} the condition that 
$\Re p(\mathbf{T})$ ($=\Re h(\mathbf{T})$) is a positive semidefinite operator for every $N$-tuple $\mathbf{T}$ of $\Lambda$-jointly nilpotent contractive $n\times n$ matrices, for all $n\in\mathbb{N}$, is fulfilled. By Theorem~\ref{thm:main}, there exists $f\in\mathcal{HA}_N^{\rm nc}(\mathcal{Y})$ such that $f_\emptyset=\frac{c_\emptyset}{2}=h_\emptyset$, $f_w=c_w=h_w$ for $w\in\Lambda\setminus\{\emptyset\}$. Then the formal power series $F(z):=(f(z)-I_\mathcal{Y})(f(z)+I_\mathcal{Y})^{-1}\in\mathcal{L(Y)}\left\langle \left\langle z_1,\ldots,z_N\right\rangle\right\rangle$ is well defined and belongs to the class $\mathcal{SA}_N^{\rm nc}(\mathcal{Y})$ (see Section~\ref{sec:ha}). Since $f(\mathbf{T})=h(\mathbf{T})=p(\mathbf{T})$  for any $N$-tuple $\mathbf{T}$ of $\Lambda$-jointly nilpotent contractive $n\times n$ matrices,  $n\in\mathbb{N}$, we get $F(\mathbf{T})=q(\mathbf{T})$ for  such a $\mathbf{T}$. Thus, the polynomial
$\sum_{w\in\Lambda}(F_w-s_w)z^w\in\mathcal{L(Y)}\left\langle z_1,\ldots,z_N\right\rangle$ vanishes on $N$-tuples of $\Lambda$-jointly nilpotent contractive $n\times n$ matrices, for every $n\in\mathbb{N}$. By Proposition~\ref{prop:vanish},
$F_w=s_w$ for all $w\in\Lambda$, i.e., $F$ solves Problem~\ref{prob:cf} for the data $s_w,\ w\in\Lambda$. 

Consider now the case where $\mathcal{U}=\mathcal{Y}$, however $s_\emptyset$ is not necessarily selfadjoint and negative semidefinite. The operator $s_\emptyset$ has a polar deconposition $s_\emptyset=UR$, where $U\in\mathcal{L(Y)}$ is unitary and $R\in\mathcal{L(Y)}$ is positive semidefinite. Suppose that $\| q(\mathbf{T})\|\leq 1$ holds for every $N$-tuple $\mathbf{T}$ of $\Lambda$-jointly nilpotent contractive $n\times n$ matrices, for all $n\in\mathbb{N}$. Then $s_\emptyset$ is a contraction, and $-I_\mathcal{Y}\leq-R\leq 0$. Define the operators $\widetilde{s_w}:=-U^*s_w,\ w\in\Lambda,$ and the polynomial $\widetilde{q}(z):=-U^*q(z)$. Clearly,  $\| \widetilde{q}(\mathbf{T})\|\leq 1$ holds for every $N$-tuple $\mathbf{T}$ of $\Lambda$-jointly nilpotent contractive $n\times n$ matrices, for all $n\in\mathbb{N}$, and $-I_\mathcal{Y}\leq\widetilde{s_\emptyset}=-R\leq 0$. By the result of the previous paragraph, there exists a solution $\widetilde{F}(z)\in\mathcal{SA}_N^{\rm nc}(\mathcal{Y})$ to Problem~\ref{prob:cf} for the data $\widetilde{s_w},\ w\in\Lambda$. Then $F(z):=-U\widetilde{F}(z)\in\mathcal{SA}_N^{\rm nc}(\mathcal{Y})$ is a solution to Problem~\ref{prob:cf} for the data $s_w,\ w\in\Lambda$.

Consider now the case where  $\mathcal{U}$ does not necessarily coincide with $\mathcal{Y}$.  Suppose that $\| q(\mathbf{T})\|\leq 1$ holds for every $N$-tuple $\mathbf{T}$ of $\Lambda$-jointly nilpotent contractive $n\times n$ matrices, for all $n\in\mathbb{N}$. Define \begin{eqnarray*}\widetilde{s_w} &:=& \left[\begin{array}{cc}
0 & 0\\
s_w & 0\end{array}\right]\in\mathcal{L(U\oplus Y)},\quad w\in\Lambda,\\
 \widetilde{q}(z) &:=& \left[\begin{array}{cc}
0 & 0\\
q(z) & 0\end{array}\right]\in\mathcal{L(U\oplus Y)}\left\langle z_1,\ldots,z_N\right\rangle.
\end{eqnarray*}
Clearly, $\| \widetilde{q}(\mathbf{T})\|\leq 1$ holds for every $N$-tuple $\mathbf{T}$ of $\Lambda$-jointly nilpotent contractive $n\times n$ matrices, for all $n\in\mathbb{N}$. By the result of the previous paragraph, there exists a solution $\widetilde{F}(z)\in\mathcal{SA}_N^{\rm nc}(\mathcal{U\oplus Y})$ to Problem~\ref{prob:cf} for the data $\widetilde{s_w},\ w\in\Lambda$. Then $F(z):=P_\mathcal{Y}\widetilde{F}(z)\big|_\mathcal{U}\in\mathcal{SA}_N^{\rm nc}(\mathcal{U,Y})$ is a solution to Problem~\ref{prob:cf} for the data $s_w,\ w\in\Lambda$.
\end{proof}
\begin{rem}
The referee suggested that the fact that the quotient of an operator algebra by an ideal
 is itself an operator algebra (see \cite{BRS}) may provide an alternate approach to proving Theorem~\ref{thm:cf}. We leave it now for a possible further exploration.
\end{rem}
Let us remark also that the examples analogous to Examples~\ref{ex:1} and \ref{ex:2} can be easily constructed for the setting of the present section, and one can ask the following questions. For which admissible sets $\Lambda\subset\mathcal{F}_N$ and data $s_w,\ w\in\Lambda$, the condition that $\| q(\mathbf{T})\|\leq 1$ for every $N$-tuple $\mathbf{T}$ of $n\times n$ matrices $\mathbf{T}\in\mathcal{C}^N\cap {\rm Nilp}_N(\Lambda_m),\ n\in\mathbb{N}$, where $m=\max_{w\in\Lambda}|w|$, is necessary for the solvability of Problem~\ref{prob:cf}?  Which admissible sets $\Lambda\subset\mathcal{F}_N$ are maximal in the sense that, for a certain choice of problem data $s_w,\ w\in\Lambda$, the condition above fails not only for $\Lambda_m\supset\Lambda$ but also for every admissible set $\tilde{\Lambda}\supset\Lambda$?
\section*{Aknowledgements}

I am thankful to Prof. Hugo Woerdeman for discussions which stimulated my interest to non-commutative interpolation problems. I express my gratitude to Prof. Victor Katsnelson for a copy of his unpublished monograph \cite{Kats} which is a good source of mathematical and historical information on classical interpolation problems and other related questions in analysis.


\begin{thebibliography}{10}

\bibitem{Ag2}
J.~Agler.
\newblock Interpolation.
\newblock Unpublished preprint.

\bibitem{Ag3}
J.~Agler.
\newblock Nevanlinna-{P}ick interpolation on {S}obolev space.
\newblock {\em Proc. Amer. Math. Soc.}, 108(2):341--351, 1990.

\bibitem{Ag1}
J.~Agler.
\newblock On the representation of certain holomorphic functions defined on a
  polydisc.
\newblock In {\em Topics in operator theory: Ernst D. Hellinger memorial
  volume}, volume~48 of {\em Oper. Theory Adv. Appl.}, pages 47--66.
  Birkh\"auser, Basel, 1990.

\bibitem{AgM1}
J.~Agler and J.~E. McCarthy.
\newblock Nevanlinna-{P}ick interpolation on the bidisk.
\newblock {\em J. Reine Angew. Math.}, 506:191--204, 1999.

\bibitem{AgM2}
J.~Agler and J.~E. McCarthy.
\newblock Complete {N}evanlinna-{P}ick kernels.
\newblock {\em J. Funct. Anal.}, 175(1):111--124, 2000.

\bibitem{AD}
L.~A. A{\u\i}zenberg and {\v{S}}.~A. Dautov.
\newblock Holomorphic functions of several complex variables with nonnegative
  real part. {T}races of holomorphic and pluriharmonic functions on the \v
  {S}ilov boundary.
\newblock {\em Mat. Sb. (N.S.)}, 99(141)(3):342--355, 479, 1976.

\bibitem{Akh}
N.~I. Akhiezer.
\newblock {\em The classical moment problem and some related questions in
  analysis}.
\newblock Translated by N. Kemmer. Hafner Publishing Co., New York, 1965.

\bibitem{AK2}
D.~Alpay and D.~S. Kalyuzhny\u{\i}-Verbovetzki\u{\i}.
\newblock Matrix-{$J$}-unitary non-commutat\-ive rational formal power series.
\newblock {\em Oper. Theory Adv. Appl.} To appear.

\bibitem{AT}
C.-G. Ambrozie and D.~Timotin.
\newblock A von {N}eumann type inequality for certain domains in
  {$\mathbf{C}^n$}.
\newblock {\em Proc. Amer. Math. Soc.}, 131(3):859--869 (electronic), 2003.

\bibitem{AP}
A.~Arias and G.~Popescu.
\newblock Noncommutative interpolation and {P}oisson transforms.
\newblock {\em Israel J. Math.}, 115:205--234, 2000.

\bibitem{Ar}
D.~Z. Arov.
\newblock Passive linear steady-state dynamical systems.
\newblock {\em Sibirsk. Mat. Zh.}, 20(2):211--228, 457, 1979.

\bibitem{Arv2}
W.~Arveson.
\newblock Subalgebras of {$C\sp{\ast} $}-algebras. {II}.
\newblock {\em Acta Math.}, 128(3-4):271--308, 1972.

\bibitem{Arv1}
W.~B. Arveson.
\newblock Subalgebras of {$C\sp{\ast} $}-algebras.
\newblock {\em Acta Math.}, 123:141--224, 1969.

\bibitem{BB2}
J.~A. Ball and V.~Bolotnikov.
\newblock On a bitangential interpolation problem for contractive-valued
  functions on the unit ball.
\newblock {\em Linear Algebra Appl.}, 353:107--147, 2002.

\bibitem{BB1}
J.~A. Ball and V.~Bolotnikov.
\newblock Realization and interpolation for {S}chur-{A}gler-class functions on
  domains with matrix polynomial defining function in {$\mathbb{C}^n$}.
\newblock {\em J. Funct. Anal.}, 213(1):45--87, 2004.

\bibitem{BGM}
J.~A. Ball, G.~Groenewald, and T.~Malakorn.
\newblock Conservative structured noncommutative multidimensional linear
  systems.
\newblock Preprint, 2003.

\bibitem{BLTT}
J.~A. Ball, W.~S. Li, D.~Timotin, and T.~T. Trent.
\newblock A commutant lifting theorem on the polydisc.
\newblock {\em Indiana Univ. Math. J.}, 48(2):653--675, 1999.

\bibitem{BT}
J.~A. Ball and T.~T. Trent.
\newblock Unitary colligations, reproducing kernel {H}ilbert spaces, and
  {N}evanlinna-{P}ick interpolation in several variables.
\newblock {\em J. Funct. Anal.}, 157(1):1--61, 1998.

\bibitem{BTV}
J.~A. Ball, T.~T. Trent, and V.~Vinnikov.
\newblock Interpolation and commutant lifting for multipliers on reproducing
  kernel {H}ilbert spaces.
\newblock In {\em Operator theory and analysis (Amsterdam, 1997)}, volume 122
  of {\em Oper. Theory Adv. Appl.}, pages 89--138. Birkh\"auser, Basel, 2001.

\bibitem{BRS}
D.~P.~Blecher, Z.-J.~Ruan, and A.~M.~Sinclair.
\newblock A characterization of operator algebras.
\newblock {\em J. Funct. Anal.}, 89(1):188--201, 1990.


\bibitem{Car1}
C.~Carath\'{e}odory.
\newblock \"{U}ber den {V}ariabilit\"{a}tsbereich der {K}oeffizienten von
  {P}otenzreihen, die gegebene {W}erte nicht annehmen.
\newblock {\em Math. Ann.}, 64:95--115, 1907.

\bibitem{Car2}
C.~Carath\'{e}odory.
\newblock \"{U}ber den {V}ariabilit\"{a}tsbereich der {F}ourier'schen
  {K}onstanten von positiven harmonischen {F}unktionen.
\newblock {\em Rend. Circ. Mat. Palermo}, 32:193--217, 1911.

\bibitem{CF}
C.~Carath\'{e}odory and L.~Fej\'{e}r.
\newblock \"{U}ber den {Z}usammenhang der {E}xtremen von harmonischen
  {F}unktionen mit ihren {K}oeffizienten und \"{u}ber den
  {P}icard--{L}andau'schen {S}atz.
\newblock {\em Rend. Circ. Mat. Palermo}, 32:218--239, 1911.

\bibitem{CW}
B.~J. Cole and J.~Wermer.
\newblock Pick interpolation, von {N}eumann inequalities, and hyperconvex sets.
\newblock In {\em Complex potential theory (Montreal, PQ, 1993)}, volume 439 of
  {\em NATO Adv. Sci. Inst. Ser. C Math. Phys. Sci.}, pages 89--129. Kluwer
  Acad. Publ., Dordrecht, 1994.

\bibitem{Co}
T.~Constantinescu.
\newblock {\em Schur parameters, factorization and dilation problems},
  volume~82 of {\em Operator Theory: Advances and Applications}.
\newblock Birkh\"auser Verlag, Basel, 1996.

\bibitem{CoJ}
T.~Constantinescu and J.~L. Johnson.
\newblock A note on noncommutative interpolation.
\newblock {\em Canad. Math. Bull.}, 46(1):59--70, 2003.

\bibitem{CS}
M.~Cotlar and C.~Sadosky.
\newblock Transference of metrics induced by unitary couplings, a {S}arason
  theorem for the bidimensional torus, and a {S}z.-{N}agy-{F}oias theorem for
  two pairs of dilations.
\newblock {\em J. Funct. Anal.}, 111(2):473--488, 1993.

\bibitem{DKh}
Sh.~A. Dautov and G.~Khuda{\u\i}berganov.
\newblock The {C}arath\'eodory-{F}ej\'er problem in higher-dimensional complex
  analysis.
\newblock {\em Sibirsk. Mat. Zh.}, 23(2):58--64, 215, 1982.

\bibitem{DP}
K.~R. Davidson and D.~R. Pitts.
\newblock Nevanlinna-{P}ick interpolation for non-commutative analytic
  {T}oeplitz algebras.
\newblock {\em Integral Equations Operator Theory}, 31(3):321--337, 1998.

\bibitem{EPP}
J.~Eschmeier, L.~Patton, and M.~Putinar.
\newblock Carath\'eodory-{F}ej\'er interpolation on polydisks.
\newblock {\em Math. Res. Lett.}, 7(1):25--34, 2000.

\bibitem{FF}
C.~Foias and A.~E. Frazho.
\newblock {\em The commutant lifting approach to interpolation problems},
  volume~44 of {\em Operator Theory: Advances and Applications}.
\newblock Birkh\"auser Verlag, Basel, 1990.

\bibitem{He}
G.~Herglotz.
\newblock \"{U}ber {P}otenzreien mit positiven, reelen {T}eil im
  {E}inheitkreis.
\newblock {\em Berichte \"{u}ber die Verhandlungen der k\"{o}niglich
  s\"{a}chsischen Gesellschaft der Wissenschaften zu Leipzig.
  Mathematische-physische Klasse}, 63:501--511, 1911.

\bibitem{K2}
D.~S. Kalyuzhniy.
\newblock Multiparametric dissipative linear stationary dynamical scattering
  systems: discrete case.
\newblock {\em J. Operator Theory}, 43(2):427--460, 2000.

\bibitem{KV}
D.~S. Kalyuzhny\u{\i}-Verbovetzki\u{\i} and V.~Vinnikov.
\newblock Non-commutative positive kernels and their matrix evaluations.
\newblock {\em Proc. Amer. Math. Soc.}
\newblock To appear.

\bibitem{Kats}
V.~Katsnelson.
\newblock Extremal problems of {G. Szeg\"{o} and M. Riesz}, factorization
  problems and other related questions in analysis. {Part I. The} scalar case.
\newblock Unpublished monograph, 1991.

\bibitem{KP}
A.~Kor\'{a}nyi and L.~Puk\'{a}nszky.
\newblock Holomorphic functions with positive real part on polycylinders.
\newblock {\em Trans. Amer. Math. Soc.}, 108:449--456, 1963.

\bibitem{KN}
M.~G. Kre{\u\i}n and A.~A. Nudel{'}man.
\newblock {\em The {M}arkov moment problem and extremal problems}.
\newblock American Mathematical Society, Providence, R.I., 1977.
\newblock Ideas and problems of P. L. \v Ceby\v sev and A. A. Markov and their
  further development, Translated from the Russian by D. Louvish, Translations
  of Mathematical Monographs, Vol. 50.

\bibitem{McC1}
S.~McCullough.
\newblock Carath\'eodory interpolation kernels.
\newblock {\em Integral Equations Operator Theory}, 15(1):43--71, 1992.

\bibitem{McC2}
S.~McCullough.
\newblock The local de {B}ranges-{R}ovnyak construction and complete
  {N}evanlinna-{P}ick kernels.
\newblock In {\em Algebraic methods in operator theory}, pages 15--24.
  Birkh\"auser Boston, Boston, MA, 1994.

\bibitem{McC3}
S.~McCullough.
\newblock Factorization of operator-valued polynomials in several non-commuting
  variables.
\newblock {\em Linear Algebra Appl.}, 326(1-3):193--203, 2001.

\bibitem{MS}
P.~Muhly and B.~Solel.
\newblock Hardy algebras, $W\sp *$-correspondences and interpolation theory.
\newblock {\em  Math. Ann.}, 330(2):353--415, 2004.

\bibitem{Ne}
M.~Neumark.
\newblock Positive definite operator functions on a commutative group.
\newblock {\em Bull. Acad. Sci. URSS S\'er. Math. [Izvestia Akad. Nauk SSSR]},
  7:237--244, 1943.

\bibitem{Pf}
A.~Pfister.
\newblock \"{U}ber das {K}oeffizientenproblem der beschr\"ankten {F}unktionen
  von zwei {V}er\"an\-derlichen.
\newblock {\em Math. Ann.}, 146:249--262, 1962.

\bibitem{P1}
G.~Popescu.
\newblock Multi-analytic operators and some factorization theorems.
\newblock {\em Indiana Univ. Math. J.}, 38(3):693--710, 1989.

\bibitem{P2}
G.~Popescu.
\newblock On intertwining dilations for sequences of noncommuting operators.
\newblock {\em J. Math. Anal. Appl.}, 167(2):382--402, 1992.

\bibitem{P3}
G.~Popescu.
\newblock Multi-analytic operators on {F}ock spaces.
\newblock {\em Math. Ann.}, 303(1):31--46, 1995.

\bibitem{P4}
G.~Popescu.
\newblock Interpolation problems in several variables.
\newblock {\em J. Math. Anal. Appl.}, 227(1):227--250, 1998.

\bibitem{Q}
P.~Quiggin.
\newblock For which reproducing kernel {H}ilbert spaces is {P}ick's theorem
  true?
\newblock {\em Integral Equations Operator Theory}, 16(2):244--266, 1993.

\bibitem{Ri}
F.~Riesz.
\newblock Sur certains syst\`{e}mes singuliers d'\'{e}quations int\'{e}grales.
\newblock {\em Annales Scientifiques de l'Ecole Normale sup\'{e}rieure
  (Paris)}, 28:33--62, 1911.

\bibitem{Ro}
L.~H. Rowen.
\newblock {\em Polynomial identities in ring theory}, volume~84 of {\em Pure
  and Applied Mathematics}.
\newblock Academic Press Inc. [Harcourt Brace Jovanovich Publishers], New York,
  1980.

\bibitem{Sakh}
L.~A. Sakhnovich.
\newblock {\em Interpolation theory and its applications}, volume 428 of {\em
  Mathematics and its Applications}.
\newblock Kluwer Academic Publishers, Dordrecht, 1997.

\bibitem{Sar2}
D.~Sarason.
\newblock On spectral sets having connected complement.
\newblock {\em Acta Sci. Math. (Szeged)}, 26:289--299, 1965.

\bibitem{Sar}
D.~Sarason.
\newblock Generalized interpolation in {$H\sp{\infty }$}.
\newblock {\em Trans. Amer. Math. Soc.}, 127:179--203, 1967.

\bibitem{Sch}
I.~Schur.
\newblock \"{U}ber {P}otenzreihen die im {I}nnern des {E}inheitskreises
  beschr\"{a}nkt sind.
\newblock {\em J. Reine Angew. Math.}, 147:205--232, 1917.

\bibitem{Sh}
B.~V. Shabat.
\newblock {\em Introduction to complex analysis. {P}art {II}}, volume 110 of
  {\em Translations of Mathematical Monographs}.
\newblock American Mathematical Society, Providence, RI, 1992.
\newblock Functions of several variables, Translated from the third (1985)
  Russian edition by J. S. Joel.

\bibitem{St}
W.~F. Stinespring.
\newblock Positive functions on {$C\sp *$}-algebras.
\newblock {\em Proc. Amer. Math. Soc.}, 6:211--216, 1955.

\bibitem{SzNF}
B.~Sz.-Nagy and C.~Foia\c{s}.
\newblock {\em Harmonic analysis of operators on {H}ilbert space}.
\newblock Translated from the French and revised. North-Holland Publishing Co.,
  Amsterdam, 1970.

\bibitem{Toep}
O.~Toeplitz.
\newblock \"{U}ber die {F}ourier'sche {E}ntwickelung positiver {F}unktionen.
\newblock {\em Rend. Circ. Mat. Palermo}, 32:191--192, 1911.

\bibitem{W}
H.~J. Woerdeman.
\newblock Positive {C}arath\'eodory interpolation on the polydisc.
\newblock {\em Integral Equations Operator Theory}, 42(2):229--242, 2002.

\end{thebibliography}

\end{document}